\documentclass[9pt, mybib]{article}

\usepackage{amssymb}
\usepackage{latexsym}
\input epsf

\def\sref#1{Section~\ref{#1}}
\def\sbref#1{Subsection~\ref{#1}}
\def\tref#1{Theorem~\ref{#1}}
\def\dref#1{Definition~\ref{#1}}

\def\rref#1{Remark~\ref{#1}}
\def\cref#1{Corollary~\ref{#1}}
\def\pref#1{Proposition~\ref{#1}}
\def\lref#1{Lemma~\ref{#1}}

\def\oref#1{Observation~\ref{#1}}
\def\fref#1{formula~(\ref{#1})}

\newtheorem{thm}{Theorem}[section] 
\newtheorem{dfn}[thm]{Definition} 
 
\newtheorem{rmk}[thm]{Remark}

\newtheorem{cor}[thm]{Corollary}
\newtheorem{prop}[thm]{Proposition} 
\newtheorem{lem}[thm]{Lemma}

\newtheorem{ex}[thm]{Example} 
\newtheorem{ob}[thm]{Observation} 
 
\newtheorem{ccl}[thm]{Conclusion} 

\newcommand{\intro}[1]
{\renewcommand{\thesection}{\fnsymbol{section}}
\setcounter{section}{-1}
\section{#1}
\renewcommand{\thesection}{\arabic{section}}
}

\def\fr#1#2{\frac{#1}{#2}}
\def\frak#1{\mathfrak{#1}}
\def\ti#1{\tilde{#1}}
\def\sq#1#2{{#1}_1, \ldots, {#1}_{#2}}
\def\dim#1{{\rm #1}\,}

\newcommand{\Pf}{{\em Proof}. }

\newcommand{\cqfd}
{%
\mbox{}%
\nolinebreak%
\hfill%
\rule{2mm}{2mm}%
\medbreak%
\par%
}
\renewcommand{\a}{{\cal A}{}} 
\renewcommand{\b}{{\cal B}{}}

\renewcommand{\d}{{\cal D}{}}  
\newcommand{\e}{{\cal E}{}}
\newcommand{\F}{{\cal F}{}}
\newcommand{\G}{{\cal G}{}}
\newcommand{\h}{{\cal H}{}} 
\newcommand{\I}{{\cal I}{}} 
\newcommand{\J}{{\cal J}{}} 
 
\newcommand{\N}{{\cal N}{}}
\newcommand{\M}{{\cal M}{}}
\newcommand{\OO}{{\cal O}{}} 
\newcommand{\p}{{\cal P}{}}

\newcommand{\R}{\mathbb R}
\newcommand{\ES}{{\cal S}{}}

\newcommand{\U}{{\cal U}{}}
\newcommand{\EK}{{\cal L}{}}
\newcommand{\KK}{{\cal K}{}} 
\newcommand{\V}{{\cal V}{}}
\newcommand{\W}{{\cal W}{}}

\newcommand{\eps}{\varepsilon }
\newcommand{\ol}{\overline}
\newcommand{\ul}{\underline}

\newcommand{\hp}{h-principle }
\newcommand{\tdco}{tangential de Rham cohomology }
\newcommand{\td}{tangential differential }
\newcommand{\s}{symplectic }

\newcommand{\lw}{leafwise }

\newcommand{\lsst}{leafwise symplectic structure }
\newcommand{\lsste}{leafwise symplectic structure}

\newcommand{\fm}{foliated manifold }
\newcommand{\nd}{neighborhood }
\newcommand{\nds}{neighborhoods }
\newcommand{\tnd}{tubular neighborhood }
\newcommand{\tnds}{tubular neighborhoods }	
\newcommand{\Iff}{if and only if }

\newcommand{\Wloge}{without loss of generality}
\newcommand{\wrt}{with respect to }
\newcommand{\rp}{respectively }

\begin{document} 				

\title{A h-principle for open relations invariant under foliated
isotopies} 
\author{M\'elanie Bertelson \thanks{Max-Planck-Institut f\"ur
Mathematik, D-53111 Bonn, Germany (bertel@mpim-bonn.mpg.de). This work
has been supported by an Alfred P.~Sloan Dissertation Fellowship.}} 
\date{\today}

\maketitle

\begin{abstract}

This paper presents a natural extension to foliated spaces of the 
following result due to Gromov~: the h-principle for open, invariant
differential relations is valid on open manifolds. The definition of
{\em openness} for foliated spaces adopted here involves a certain type of
Morse functions. Consequences concerning the problem of existence of
regular Poisson structures, the original motivation for this work, are
presented. 

\end{abstract}

\intro{Introduction}

Gromov proved in \cite{Gv-69} the following theorem.  
\begin{thm}\label{Bach}
On an open manifold, the \hp for open, invariant differential 
relations is valid. 
\end{thm}
Heuristically, a {\em differential relation} on a manifold $M$ is a 
differential constraint on the sections of a certain bundle $\pi : E
\to M$. More precisely, it is a subset $\Omega$ of some jet bundle   
$J^k(E)$ of local sections of $E$. A section of $E$ whose $k$-jet 
extension is entirely contained in the relation satisfies the
constraint; it  
is called a {\em solution} of the differential relation $\Omega$.
A relation $\Omega$ is said to be {\em open} (\rp {\em invariant})
when $\Omega$ is an open subset of $J^k(E)$ (\rp when isotopies of $M$
can be lifted to isotopies of $J^k(E)$ that preserve
$\Omega$). Letting $Sol(\Omega)$ (\rp $\Gamma(\Omega)$) denote the set   
of solutions (\rp sections) of $\Omega$, endowed with the weak $C^k$
(\rp $C^0$) topology, the relation $\Omega$ is said {\em to satisfy
the h-principle} if the $k$-jet extension map     
$$
j^k : Sol(\Omega) \to \Gamma(\Omega) : f \mapsto j^kf\;
$$
is a {\em weak homotopy equivalence}. This means that the map $j^k$ 
induces a bijection inbetween arcwise connected components of
$Sol(\Omega)$ and $\Gamma(\Omega)$, and isomorphisms of the homotopy
groups for the various components.

\tref{Bach} has a lot of important corollaries concerning the existence 
and classification problems of various types of objects in differential 
topology and geometry (cf.~\cite{Gv-book}). For instance, applying 
\tref{Bach} to the relation $\Omega^{\cal S} = \{\, j^1\alpha(x) \in 
J^1(T^*M); d\alpha(x)$ is nondegenerate $\}$ yields the following 
result. 

\begin{cor}\label{Massenet} Let $M$ be an open manifold. The inclusion
of the space of exact symplectic forms on $M$ into the space of
nondegenerate $2$-forms is a weak homotopy equivalence. 
\end{cor}
In particular, on an open manifold, existence of a symplectic
structure depends only on existence of a nondegenerate $2$-form, a
problem that belongs to obstruction theory. \\ 

The hypothesis that $M$ is open is crucial. Indeed, already for the
symplectic relation $\Omega^{\cal S}$, the h-principle is far from
being valid on a closed manifold. In addition to a nondegenerate
$2$-form, a closed symplectic manifold admits a de~Rham class in
$H^2(M)$ whose top exterior power does not vanish. More subtle 
conditions, involving Seiberg--Witten invariants, have been discovered
by Taubes (cf.~\cite{CT}). Furthermore, even when the manifold admits
a symplectic structure, not any nondegenerate $2$-form may be deformed
into a symplectic form. In general, if every open, invariant
differential relation defined on a manifold $M$ satisfies the 
h-principle, then $M$ must be open, as such a relation can be
constructed that admits sections and whose solutions are functions
without local maxima.\\   

Motivated by the problem of existence of \lw \s structures on foliated
spaces (cf.~\cite{B, B-f}), we searched for a generalization of
\tref{Bach} to {\em foliated invariant} differential relations, that
is, differential relations that are invariant under isotopies that
preserve a certain foliation on the manifold. This requires
finding a good notion of ``openness'' for foliated spaces. It is
important to observe that one may not, in general, impose on the
solutions constructed in the proof of \tref{Bach} to be nicely behaved
at infinity. In contrast, a solution of a differential relation (a
symplectic structure for instance) on a nonclosed leaf of a foliation,
that is the restriction of a global solution (a leafwise symplectic
structure) is most likely very constrained at infinity, partly due to
recurrence phenomena, partly due to the influence of neighboring
leaves. The foliated case lies, in some sense, midway between the open
case and the closed case. Examples of foliations have been exhibited
in \cite{B, B-f} that seem from some point of view quite open, that
support \lw nondegenerate $2$-forms, but do not admit any \lsste. 
One should therefore avoid to include them in the class of open
foliations. \\ 

On the other hand, open manifolds are characterized by the existence 
of a positive, proper Morse function, without any local maximum. 
The proof of \tref{Bach} suggests to base a definition of openness 
for foliated manifolds on that characterization. This justifies the 
following definition.  

\begin{dfn}\label{Mahler} A foliated manifold $(M,\F)$ is said to be
{\em open} if there exists a smooth function $f: M \to [0,\infty)$ that 
has the following properties~: 
\begin{enumerate}
\item[a)] $f$ is proper,
\item[b)] $f$ has no leafwise local maxima, 
\item[c)] $f$ is $\F$-generic (cf.~\dref{cello} below).
\end{enumerate}
\end{dfn}
With this definition of openness, the following result holds~:

\begin{thm}\label{ht} On an open foliated manifold, any open, foliated
invariant differential relation satisfies the parametric h-principle. 
\end{thm}

\begin{cor}\label{Stravinsky} Let $(M,\F)$ be an open foliated manifold. 
Any \lw nondegenerate $2$-form is homotopic, in the class of \lw 
nondegenerate $2$-forms, to a leafwise symplectic form. 
\end{cor}

The leaves of an open foliation are necessarily open manifolds, but
this condition is not sufficient. One might interpret \dref{Mahler} as
that of {\em uniform openness} of the leaves. It can be checked
directly (i.e.~without quoting \cref{Stravinsky}) that the foliated
manifolds introduced in \cite{B, B-f} do not support any function $f$
satisfying $a)$ and~$b)$.\\ 

The proof of \tref{ht} involves consideration of the trajectories of a {\em \lw 
gradient vector field} for $f$, as did the proof of \tref{Bach}. 
There are new difficulties. First, the {\em leafwise critical 
points} of $f$ are not isolated but come in families. Thus they cannot
be handled one at a time (as they are in the nonfoliated
case). Secondly, leafwise critical points may be degenerate, 
even generically. The set of trajectories converging to a degenerate
critical point is not in general well understood. To overcome the
latter difficulty one needs to carry out a careful construction of a
Riemannian metric for which the trajectories of the associated \lw
gradient vector field are somewhat controlled near the singular 
locus.\\

This paper is organized as follows. \sref{hpbg} is a brief
introduction to the theory of h-principles (see \cite{Gv-book}). 
\sref{outline} presents an outline of the proof of \tref{ht}, and
should serve as a reading-guide for the remainder of the
text. \sref{gen} is devoted to making precise the term {\em 
$\F$-generic} referred to in \dref{Mahler}, and to introducing the
notion of {\em strong $\F$-genericity} needed later. \sref{lwgr}
exhibits some properties of {\em leafwise gradient vector fields}
(\dref{lwgrdef}) of strongly $\F$-generic functions. \sref{metric}
describes the construction of a nice Riemannian metric associated to a
strongly $\F$-generic function. The proof of \tref{ht} is completed in
\sref{proof}. Some examples of open foliated manifolds are presented 
in \sref{ex}. Finally, \sref{Poisson} is concerned with
\cref{Stravinsky}.\\ 

{\bf Acknowledgements.}
I wish to thank Emmanuel Giroux, who suggested to generalize Gromov's
result, Alan Weinstein, whose insightful advice has helped this work
reach its final shape, and Yasha Eliashberg for many very stimulating
conversations.  

\section{H-principles}\label{hpbg}

We state here the definitions and results of the theory of
h-principles (cf.~\cite{Gv-book}) that will be needed in the remainder
of the text. We have followed Emmanuel Giroux's beautiful (as yet
unpublished) lecture notes~\cite{Giroux}. The proofs are reproduced in
\cite{B}.     

\subsection{Differential relations and h-principles}\label{Gounod}

Consider a locally trivial fibration $\pi : E \to M$ with fiber a
manifold $F$. Let $E^k$ denote the set of $k$-jets of local sections
of $\pi : E \to M$. If $f:U\subset M\to E$ is a local section defined 
on an open subset $U$ of $M$, the $k$-jet of $f$ at $x\in U$ is 
denoted by $j^{k}f(x)$. The set $E^{k}$, endowed with the natural 
projection $E^k \to M : j^{k}f(x) \mapsto x$, is a locally trivial 
fibration. For $k\leq r\leq \infty$, a $C^{r}$ local section $f$ of
$E$, defined on an open subset $U$, induces a $C^{r-k}$ local section
$j^kf : U \to E^k : x \mapsto j^kf(x)$ of $E^k$, called the {\em $k$-jet 
extension of~$f$}. 

\begin{dfn}
A subset $\Omega$ of $E^k$ is called a {\sl differential relation of
order $k$}. It is said to be {\em open} if it is an open subset.
A section of $E^k$ whose values are in $\Omega$ is called a {\sl
section of $\Omega$}. A $C^k$ section of $E$ whose $k$-jet extension is a
section of $\Omega$ is called a {\sl solution of $\Omega$}. A section
of $E^k$ that coincides with the $k$-jet extension of some section of
$E$ is said to be {\em holonomic}.  
\end{dfn}
Observe that solutions of $\Omega$ and holonomic sections of $\Omega$
are in one to one correspondence.\\

Let $\Omega$ be a differential relation on the manifold $M$. We will
be considering families of sections of $\Omega$ parameterized by cubes
$S = [0,1]^p$,$p \geq 1$. Those are defined to be continuous maps $f :
M \times S \to E^k$, such that for each $s$ in $S$, the partial map
$f_s : M \to E^k$ (obtained by restricting the map $f$ to $M\times
\{s\}$) is a section of $\Omega$. A {\em homotopy of sections of
$\Omega$} is a family parameterized by $S=[0,1]$. Concerning local
sections of $\Omega$, we need to introduce some terminology. Let
$(A,A')$ denote a nested pair of compact subsets of $M$. The word
``nested'' indicates that $A' \subset A$. 

\begin{enumerate}
\item[-] A {\em family of sections of $\Omega$ defined near $A$} is a
family defined on $U \times S$ for some \nd $U$ of~$A$.
\item[-] Two families of sections $f_s$  and $g_s$ are said to {\em
coincide near $A$} if there exists a neighborhood $U$ of $A$ on
which both $f_s$ and $g_s$ are defined, and for which $f_s|_U = g_s|_U$
for all $s$ in $S$. 
\item[-] A family of sections $g_s$ defined near $A$ is said to {\em
extend} another family $f_s$ defined near $A'$ if $g_s$ and $f_s$
coincide near $A'$.
\item[-] Two families of sections $f_s$ and $g_s$ defined near $A$ are
said to be {\em homotopic} if there exists a third family
$h_{s,t}$ defined near $A$ with $(s,t)$ in $S\times [0,1]$, such that
$h_{s,0}$ coincides with $f_s$ near $A$ and $h_{s,1}$ coincides with
$g_s$ near~$A$.  
\item[-] A homotopy $h_{s,t}$ is said to be {\em stationary near $A$ for 
$s$ in $S'\subset S$} if there exists a neighborhood $U$ of $A$ for which 
$h_{s,t} = h_{s,0}$ on $U$, for $t$ in $[0,1]$, and for $s$ in~$S'$.
\end{enumerate}
\noindent
Similar definitions apply to sections of $E$ and to solutions of 
$\Omega$, with one restriction~: families $f : M \times S \to E$ of
$C^k$ sections of $E$ are required to be $C^k$-continuous, that is, their
$k$-jet extension $j^kf : M \times S \to E^k$ is required to be a
continuous map. To obtain solutions that are of smoothness class
$C^r$, with $k\leq r \leq \infty$, one should everywhere consider
$C^r$ sections of $E$ and $C^{k-r}$ sections of $E^k$ only. Also, when
the relations considered is open, any (continuous) family of $C^r$
sections of $E$ can be approximated by a smooth family of $C^r$ sections,
that is a smooth map $M \times S \to E$.\\  

In the sequel, $S'$ will always denote a subset of $S = [0,1]^p$
consisting of the union of some of its faces.  

\begin{dfn}(The parametric h-principle.)
The relation $\Omega\subset E^k$ is said to {\em satisfy the parametric 
h-principle on $M$ (near a subset $A$ of $M$)} if any family of sections 
of $\Omega$ (defined near $A$) is homotopic, among families of sections 
of $\Omega$ (defined near $A$), to a family of holonomic sections of $\Omega$. 
If the sections are already holonomic on $M$ (or near $A$) for $s$ in 
$S'$, we may assume them to remain holonomic during the homotopy. 
(Equivalently, the homotopy may be assumed to be stationary for s in~$S'$).  
\end{dfn}

\begin{rmk}\label{clavecin}{\rm Observe that if a differential relation
$\Omega$ satisfies the parametric h-principle, then its solutions satisfy 
some kind of uniqueness property. Indeed, let $f_0$ and $f_1$ be two
solutions of $\Omega$ whose $k$-jet extensions are homotopic among
sections of $\Omega$. Let $f^1_t$ be such a homotopy. Then the family
$f^1_t$ is homotopic to a family $f_t$ of holonomic sections of
$\Omega$ via a homotopy that is stationary for $t \in \{0,1\}$. In
particular, the two solutions $f_0$ and $f_1$ are homotopic among
solutions of~$\Omega$.  
}\end{rmk}

\begin{rmk}{\rm
The relation $\Omega$ satisfies the parametric \hp on $M$ \Iff the map
$j^k : Sol(\Omega) \to \Gamma (\Omega)$ is a weak homotopy
equivalence, where $Sol(\Omega)$ is the set of global solutions of
$\Omega$ and where $\Gamma (\Omega)$ is the set of global sections of
$\Omega$ (cf.~\cite{Gv-book}~$(C)$ p.~16). 
}\end{rmk}

\subsection{Invariant relations}\label{inv}

\begin{dfn}\label{Jongen}
An {\em isotopy of the manifold $M$} is a family $\varphi_t$,
$t$ in $[0,1]$, of diffeomorphisms of $M$ such that the map $\varphi :
M \times [0,1] \to M : (x,t) \mapsto \varphi_t(x)$ is smooth and
$\varphi_0=\dim{Id}_{M}$. Consider a foliation $\F$ on $M$. A {\em
foliated isotopy of $(M,\F)$} is an isotopy $\varphi_t$ of $M$ that
preserves the foliation $\F$, that is, $(\varphi_t)_*(T\F) = T\F$ for
all $t$ in $[0,1]$. Finally, two sets $A$ and $A'$, with $A' \subset
A$ are said to be {\em isotopic} (\rp {\em foliated isotopic)} if for
every \nd $U$ of $A'$, there exists an isotopy (\rp a foliated
isotopy) $\varphi_t^U$ of $M$ such that 
\begin{enumerate}
\item[-] $\varphi^{U}_t$ is stationary near $A'$, 
\item[-] $\varphi^{U}_t(A) \subset A$ for all $t$, 
\item[-] $\varphi^{U}_1(A)\subset U$. 
\end{enumerate}
\end{dfn}

\begin{dfn}(Invariant relations.) Let $\pi : E \to M$ be a locally
trivial fibration. An isotopy $\varphi_{t}$ of $M$ is said to {\em
operate on local sections of $E$} if we are given an isotopy
$\ol{\varphi}_{t}$ of $E$ covering $\varphi_{t}$ (i.e.~$\pi \circ
\ol{\varphi}_{t} = \varphi_{t} \circ \pi$) satisfying the following
property~:   
\begin{enumerate}
\item[] if $\varphi_{t}$ coincides with $\varphi_{t_{0}}$ on an open
        subset $U$ for all $t\geq t_{0}$, then $\ol{\varphi}_{t}$
        coincides with $\ol{\varphi}_{t_{0}}$ on $\pi^{-1}(U)$ for all
        $t\geq t_{0}$.    
\end{enumerate}
In this situation, if $f : U \to E|_{U}$ is a local section of $E$ defined
on some open subset $U$ of $M$, then $\varphi_{t}\cdot f = \ol 
\varphi_{t}^{-1}\circ f\circ\varphi_{t}$ is a local section of $E$
defined  on $\varphi_t^{-1}(U)$. This operation on local sections induces
an operation on $E^{k}$ as follows~: $\varphi_{t} \cdot (j^{k}f(x)) = j^{k}
(\varphi_{t}\cdot f) (\varphi_{t}^{-1}(x))$.
\end{dfn}

\begin{dfn}The differential relation $\Omega\subset E^{k}$ is said to be 
{\em invariant under the isotopy $\varphi_{t}$} if $\varphi_{t}$ 
operates on local sections of $E$, and if $\Omega$ is invariant
under the induced operation on $E^{k}$ (i.e.~$\varphi_{t}\cdot\Omega =
\Omega$). The relation $\Omega$ is said to be {\em invariant} (\rp
{\em foliated invariant}) if it is invariant under all isotopies of
$M$ (\rp all foliated isotopies of $(M,\F)$). 
\end{dfn}

\subsection{Local h-principle, h-principle for extensions}

\begin{dfn}\label{localhp}(The local parametric h-principle.) A
differential relation is said to satisfy the {\em local parametric
h-principle} if it satisfies the parametric h-principle near any
point.  
\end{dfn}

\begin{prop}[\cite{Gv-book} $B_2$ p.~37]\label{Schumann} Any open
differential relation satisfies the local parametric h-prin\-ci\-ple.  
\end{prop}

\begin{dfn}(The parametric h-principle for extensions.)
Let $(A,A')$ be a nested pair of compact subsets of $M$. A
differential relation $\Omega$ defined on $M$ is said to satisfy the
{\em parametric 
h-principle for extensions of solutions form $A'$ to $A$}, or {\em on
the pair $(A,A')$}, if any family $f_{s}$ of sections of $\Omega$
defined near $A$ and holonomic near $A'$ is homotopic to a family of
holonomic sections, through a homotopy that is stationary near
$A'$. Moreover, if the sections $f_{s}$ are already holonomic near $A$
for $s$ in $S'$, then the homotopy may be chosen to be stationary near
$A$ for $s$ in~$S'$.    
\end{dfn}

\begin{lem}[\cite{Gv-book}]\label{violon}
Let $\Omega$ be a differential relation on the manifold $M$, and let
$(A,A')$ be a nested pair of isotopic compact subsets of $M$. If
$\Omega$ is invariant under all isotopies $\varphi^{U}_t$, where $U$
runs through the set of \nds of $A'$, then the parametric 
h-principle for extensions holds on the pair~$(A,A')$.    
\end{lem}

\begin{prop}[\cite{Gv-book} $(A')$ p.~40]\label{exhaustion}
Let $\Omega$ be a differential relation on $M$, and let $M=\cup_{i\geq
0} K_{i}$ be an exhaustion of $M$ by compact subsets (i.e.~$K_i \subset 
K_{i+1}$ and $M = \cup_{i\geq 0} K_i$). Suppose that $\Omega$ satisfies 
the parametric \hp on $K_0$, as well as the parametric \hp for
extensions on all pairs $(K_{i+1},K_i)$. Then $\Omega$ satisfies the
h-principle on~$M$.   
\end{prop}

Consider an open differential relation $\Omega$ on a manifold $M$
endowed with a foliation $\F$. 

\begin{thm}[\cite{Gv-book} $(B_1)$ p.~41, $(C_1)$ p.~42, $(C'_3)$
p.~43]\label{Chopin}  
Let $(A,A')$ be a ne\-sted pair of compact subsets of $M$ such that
the compact $C=\ol{A-A'}$ is contained in an embedded submanifold
$M_0$ of $M$, of codimension at least one. If the relation $\Omega$ is
invariant, then it satisfies the h-principle for extensions
on~$(A,A')$. If the relation $\Omega$ is only foliated invariant, but
the submanifold $M_0$ intersects the foliation $\F$ transversely, then
$\Omega$ satisfies the h-principle for extensions on~$(A,A')$ as well. 
\end{thm}

\begin{rmk}\label{orgue}{\rm 
The hypothesis of the previous theorem can be weakened (without affecting 
the conclusion) as follows~: $C$ consists of a finite union of compact 
sets, each contained in an embedded submanifold of codimension at
least one (transverse to $\F$) (cf.~\cite{Gv-book} $(A')$ p.~40).  
}\end{rmk}

\section{Outline of the proof}\label{outline}  

We begin this section with a rough outline of the proof of
\tref{Bach}. Let $\Omega\subset E^k$ be an open, invariant
differential relation on a manifold $M$. Fixing a section $g$  
of $\Omega$, the procedure to deform $g$ into a holonomic section of 
$\Omega$ is sketched below (the case of a family of sections is
treated similarly). The hypotheses on $\Omega$ imply the 
following two facts.

\begin{enumerate}
\item[(1)] Since the relation $\Omega$ is open, it satisfies the local 
           parametric h-principle (cf.~\pref{Schumann}). Thus, for any 
           given point $x$ in $M$, the section $\beta$ can be deformed 
           into a holonomic sections on a sufficiently small \nd of $x$.    
\item[(2)] Since $\Omega$ is invariant, it satisfies the \hp for extensions 
           on any pair $(A,A')$ of isotopic compact subsets 
           (cf.~\lref{violon}). Thus, if $\beta$ can be deformed    
           into a holonomic section near $A'$, it can also be deformed 
           into a holonomic section near $A$.     
\end{enumerate}

Starting from a holonomic section $\omega$ of $\Omega$ defined on a
\nd $U$ of a point $x$ and homotopic to $\beta|_U$, one tries to
extend $\omega$ as far as possible. It can certainly be extended to
any ``large'' open ball containing $x$ (cf.~(2)). But when trying to
go further, one has to understand how to deal with the topology of
$M$. A good grasp on the latter is provided  by a proper, positive,
Morse function on $M$, that is, a proper map $f : M \to [0,\infty)$
whose singular points are nondegenerate and lie on distinct
levels. The term {\em proper} indicates that for any pair $a<b$ of real
numbers, the set $f^{-1}([a,b])$ is compact. 

Let $a < b$ be two noncritical values of $f$, and let $\omega$ be a
holonomic section of $\Omega$, defined on a \nd of $f^{-1}([0,a])$ and
homotopic to $\beta$. Provided $[a,b]$ does not contain any critical
value of $f$, the set $f^{-1}([0,b])$ is isotopic to $f^{-1}([0,a])$
(\cite{Milnor}). The holonomic section $\omega$ can therefore be
extended to a \nd of $f^{-1}([0,b])$ (cf.~$(2)$). If $[a,b]$
contains a critical value of $f$ corresponding to a critical point
$x$, the set $f^{-1}([0,b])$ is obtained from $f^{-1}([0,a])$ by
gluing a closed disk $D^k$ along its boundary, in the sense that
$f^{-1}([0,b])$ is isotopic to $f^{-1}([0,a]) \cup_{\partial D^k}
D^k$~(\cite{Milnor}). The dimension $k$ of the disk coincides with the
number of negative eigenvalues of the Hessian of $f$ at $x$. Thus, to
extend $\omega$ through $x$, one needs to be able to extend a
holonomic section of $\Omega$, homotopic to $\beta$, defined on a \nd
of the boundary of an embedded disk $D^k$, to a holonomic section of
$\Omega$, homotopic to $\beta$, defined on a \nd of the entire 
disk. (In other words, one needs the relation $\Omega$ to satisfy the
parametric \hp for extensions on the pair $(D^k,\partial D^k)$). This
can be done as long as $k < {\rm dim }\,M$ (cf.~\tref{Chopin}), and
constitutes the {\it key step} of the proof of \tref{Bach}.   

The restriction $k < {\rm dim }\,M$ explains the dichotomy between
closed and open manifolds. Indeed, if $M$ is open, we may assume that
$f$ has no local maxima, or equivalently, that the disks we glue are
never of maximal dimension. Beginning with a holonomic section
$\omega$ of $\Omega$ defined near $f^{-1}(0)$ and homotopic to
$\beta$, we can therefore extend it to $f^{-1}([0,b])$ for larger and
larger values of $b$, eventually obtaining a global holonomic section
of $\Omega$ (cf.~\pref{exhaustion}). \\  

Introducing the gradient $\nabla f$ of $f$ \wrt some Riemannian metric, 
one observes that its flow yields an isotopy between $f^{-1}([0,b])$ and 
$f^{-1}([0,a]) \cup {}_{\partial D^k} D^k$. Moreover, the embedded disk 
$D^k$ may be thought of as the set of points in $f^{-1}([a,b])$ lying 
on trajectories $\theta(t)$ ``converging'' to $x$, in the sense that 
$\lim_{t\to +\infty} \theta(t) = x$. It is useful to adopt this point of 
view when dealing with foliated manifolds. \\

Consider now the corresponding problem in the foliated case. Let $\Omega$ 
be an open, foliated invariant differential relation defined on an open 
\fm $(M,\F)$. Let $f : M \to [0,\infty)$ satisfy the hypotheses
$a)$, $b)$ and $c)$ of \dref{Mahler}. We may take advantage of what is
already known. In particular, the facts $(1)$ and $(2)$ are still true
here, provided we restrict ourselves to foliated isotopies. The {\it
key step} is also valid in this context, provided the embedded disk  
$D^k$, does not only have codimension at least one, but intersects
$\F$ transversely as well (cf.~\tref{Chopin}). The extra difficulties
one faces come from the restriction to {\em foliated}
isotopies. Indeed, up to foliated isotopy, the passage from
$f^{-1}([0,b])$ to $f^{-1}([0,b'])$ does not correspond anymore to
gluing some disk (nor even some family 
of disks). It corresponds instead to gluing the ``skeleton''
consisting of the set of bounded trajectories of the {\em leafwise
gradient vector field of $f$} associated to some Riemannian metric $g$
(cf.~\dref{lwgrdef} and \lref{alto}). The skeleton can be very 
complicated due to two phenomena~: 

\begin{enumerate}
\item[i)] the {\em \lw critical points of $f$} (\dref{lwcr}) can be \lw
degenerate (\rref{Mendelshon}),    
\item[ii)] the \lw critical points of $f$ are not isolated in $M$
(\rref{Donizetti}).     
\end{enumerate} 

Problem i) makes it hard to describe the topological structure of the
skeleton already locally, near the foliated singular locus $\Sigma_f$
(\dref{lwcr}). Fortunately, to apply the {\it key step}, we do not
need to know the exact topological type of the skeleton, it 
suffices to know that it is contained in a finite union of compact
subsets of submanifolds of codimension at least one, intersecting $\F$
transversely (cf.~\rref{orgue}). When the latter holds in a \nd of
$\Sigma_f$, the metric $g$ is said to be {\em nice}. The construction
of a nice metric is detailed in \sref{metric}. As it requires the
foliated singular locus of $f$ to be stratified according to
Thom--Boardman (\sbref{TBs}), the function $f$ will be assumed to be
{\em strongly $\F$-generic} (\dref{strongcello}). As noticed in
\rref{Verdi}, this extra assumption is not restrictive. The
construction of the metric is done by successive extensions from a \nd 
of one stratum to the next (\wrt some natural order on the set of
strata).   

Problem ii), on the other hand, makes it hard to have a grasp on the
global structure of the skeleton due to the presence of trajectories
$\theta(t)$ for which both $\lim_{t \to +\infty} \theta(t)$ and $\lim_{t
\to -\infty} \theta(t)$ are in $f^{-1}([a,b])$. For such a trajectory, 
the structure of the skeleton near $\lim_{t \to -\infty} \theta(t)$ is 
quite complicated, more so that the stable and unstable manifolds of
distinct critical points may not be assumed to intersect transversely. 
Fortunately, this difficulty vanishes if we approximate $f$ by a function 
whose \lw critical points are isolated \wrt the leaf topology 
(cf.~\pref{isolate}), as it allows one to ``slice'' $M$ sufficiently finely 
to (more or less) avoid having trajectories $\theta(t)$ for which both 
$\lim_{t \to +\infty} \theta(t)$ and $\lim_{t \to -\infty} \theta(t)$ lie 
in the same slice (cf.~\sref{proof}).    

\section{Genericity}\label{gen}

This section is organized as follows. The first subsection shows how
to adapt Boardman's construction of a natural stratification of the 
singular locus of a smooth map (\cite{Boardman}) to the foliated
case. We use Mather's description of the Thom--Boardman stratification
\cite{JM-73}. The second subsection shows that generically, the \lw
critical points of a real-valued function are isolated \wrt the leaf
topology. This is a consequence of a result mentioned in \cite{JMt}
(6.1  p.~29) and proved along a scheme that appears in \cite{Igusa}.
The third subsection defines (strongly) $\F$-generic functions and
exhibits some of their properties. We begin with recalling the
statement of Thom's transversality theorem which plays a crucial role
throughout this section, and with defining the terms {\em \lw critical
point} and {\em foliated singular locus} used in the previous section. \\    
   
Given two manifolds $M$ and $N$, the set $C^\infty(M,N)$ is endowed hereafter 
with the fine (or Whitney) $C^\infty$ topology. A subset of a topological 
space is {\em residual } if it is a countable intersection of dense open 
subsets. It is a classical result (e.g.~\cite{GG}) that a residual subset 
of $C^\infty(M,N)$ is dense. Notice that a countable intersection of residual 
sets is still residual (while an intersection of dense sets is not dense in
general). A condition on smooth functions $f$ in $C^\infty(M,N)$ is said 
to be {\em generic} if it is satisfied by all functions in a residual subset 
of $C^\infty(M,N)$. Recall that {\em a smooth map $f : M \to P$ intersects
an embedded submanifold $W$ of $P$ transversely at a point $x$ in $M$}
if and only if either $f(x) \notin W$, or $f(x)\in W$~and 
\begin{equation}\label{tr}
f_{*_x}T_xM + T_{f(x)}W = T_{f(x)}P\;,
\end{equation}
and that {\em $f$ intersects $W$ transversely} if (\ref{tr}) holds for 
all $x$ in $f^{-1}(W)$. In that situation, the set $f^{-1}(W)$ is an embedded 
submanifold of $M$. If the map $f$ is the inclusion of a submanifold $W'$, we 
say that {\em the submanifolds $W$ and $W'$ intersect transversely}. Finally, 
the set of $k$-jets of local maps $U\subset M \to N$ is denoted 
by~$J^k(M,N)$.

\begin{thm}\label{Thom}(Thom Transversality Theorem.) Let $W$ be an
embedded submanifold of $J^k(M,N)$. The set of smooth maps $f$ whose
$k$-jet extension 
$$                    
j^kf : M \to J^k(M,N)
$$
is transverse to $W$ is residual in~$C^\infty(M,N)$.
\end{thm}
A proof of the previous theorem can be found in~\cite{GG}.    

\begin{dfn}\label{lwcr} Let $(M,\F)$ be a foliated manifold, and let
$f:M\to\R$ be a smooth function. 
\begin{enumerate}
\item[-] A point $x$ at which $d(f|_{F_x})$ vanishes is
         called a {\em \lw critical point of~$f$}. 
\item[-] The {\em foliated singular locus of $f$ (\wrt the foliation $\F$)} 
         is the set of \lw critical points of $f$, it is denoted by~$\Sigma_f$.
\end{enumerate}
\end{dfn}

\begin{rmk}\label{Donizetti}{\rm An ordinary critical point of $f$ is of
course a \lw critical point as well, but \lw critical points persist in 
nearby leaves, so that they typically come in $q$-parameter families, where 
$q$ is the codimension of~$\F$.  
}\end{rmk}

\subsection{The Thom-Boardman stratification}\label{TBs}

Let $\e_n$ denote the set of germs of smooth maps $\R^n\to \R$ at
the origin. The set $\e_n$ is a local ring, whose maximal ideal, denoted 
by $\frak{m}_n$, is the set of germs vanishing at $0$. Given an ideal
$\I$ of the ring $\e_n$, we use the notation $\I^k$ for the product
$\I \cdot \ldots \cdot \I$ of $k$ copies of~$\I$. Finally, the symbol
$D_i$ stands for the derivation $\frac{\partial}{\partial x_i}$, where
$x_1, \ldots , x_n$ are the standard coordinates on $\R^n$.

\begin{dfn}Let $\I$ be a finitely generated ideal of~$\e_n$.
\begin{enumerate}
\item[i)] The {\em rank of $\,\I$}, denoted by ${\rm rk }(\I)$, is the
          dimension of $(\I + \frak{m}^2_n)/\frak{m}^2_n$ as a real
          vector space. 
\item[ii)] The ideal generated by $\I$ and the set $\Gamma_r(\I)$ of
           $r\times r$ minors of the Jacobian matrix $(D_if_j)_{1\leq
           i\leq n, 1\leq j\leq a}$, where $f_1,\ldots,f_a$ is a set of
           generators for $\I$, is denoted by $\Delta_r(\I)$ and called a
           {\em Jacobian extension of $\I$}. Notice that the ideal
           $\Delta_r(\I)$ (unlike the set $\Gamma_r(\I)$) does not
           depend on the choice of generators for~$\I$. 
\item[iii)] The ideal $\Delta_{r+1}(\I)$ with  $r=rk(\I)$ is denoted
            by $\delta(\I)$. Notice that when $\I$ is proper, the ideal
            $\delta(\I)$ is proper as~well.
\item[iv)] The {\em Boardman symbol of $\I$} is the infinite sequence
           $I(\I) = (n - r_1, n - r_2, \ldots, n - r_\ell,\ldots)$,
           where $r_\ell = {\rm rk }(\delta^{\ell-1}(\I))$.           
\end{enumerate}
\end{dfn}

Let $J^k(n,p)$ denote the set of $k$-jets of maps $(\R^n,0) \to
(\R^p,0)$ at the origin. Given a jet $z=j^kf(0)$ in $J^k(n,p)$,
represented by a map $f$, consider the ideal
$$
\I(z) \stackrel{\dim{def.}}{=} (f) + \frak{m}^{k+1}_n\;,
$$ 
where $(f)$ denotes the ideal generated by the germs of the components
$f_1,\ldots,f_p$ of $f$ at $0$. The Boardman symbol of the jet $z$,
denoted by $I(z)$, is defined to be the Boardman symbol of the ideal 
$\I(z)$ truncated at order $k$ (notice that $I(\I(z))$ is of the type
$(i_1,\ldots, i_k, 0,\ldots, 0, \ldots)$). Given a sequence $I$ of $k$
nonnegative integers, let $\Sigma^I \subset J^k(n,p)$ denote the
set of $k$-jets whose Boardman symbol is~$I$.    

\begin{prop}[\cite{Boardman}] Let $I=(\sq i k)$. The set $\Sigma^I$ is
nonempty \Iff $n\geq i_1\geq i_2\geq \ldots \geq i_k$, and either $i_1> 
n-p$, or $i_1= n-p$ and $i_1=i_2=\ldots=i_k$. Moreover $\Sigma^I$ is an 
embedded submanifold of $J^k(n,p)$ whose codimension is given by the formula 
\begin{equation}\label{cod}
{\rm cod }\,\Sigma^I = (p-n+i_1) \mu_I - (i_1-i_2) \mu_{sI} - \ldots - 
(i_{k-1} - i_{k})\mu_{s^{k-1}I} \;,
\end{equation}
where $s^jI$ denotes the sequence $(i_{j+1}, \ldots, i_k)$, and where
$\mu_{(\sq i \ell)}$ is the number of sequences
$(j_1,j_2,\ldots,j_\ell)$ of integers satisfying    
\begin{enumerate}
\item[-] $j_1\geq j_2\geq \ldots \geq j_\ell$,
\item[-] $i_r\geq j_r\geq 0$ for all $r$, and $j_1>0$.
\end{enumerate}
\end{prop}

Given two manifolds $N$ and $P$, of respective dimensions $n$ and $p$,
let $J^k(N,P)$ denote the manifold of $k$-jets of local functions on
$N$ with values in $P$. The manifold $J^k(N,P)$ is a bundle over $N
\times P$ whose fibers are diffeomorphic to $J^k(n,p)$, and whose
structure group is the product of the group of invertible jets in
$J^k(n,n)$ with the group of invertible jets in $J^k(p,p)$. For any
sequence $I=(\sq i k)$ of nonnegative integers, let $\Sigma^I \subset
J^k(N,P)$ denote the subset 
$$
\Sigma^I = \bigcup_{(x,y)\in N\times P} \Sigma^I_{x,y}\;,
$$ 
where $\Sigma^I_{x,y}$ corresponds to $\Sigma^I$ under the identification 
of the fiber of $J^k(N,P)$ at $(x,y)$ with the manifold $J^k(n,p)$. 
Since $\Sigma^I \subset J^k(n,p)$ is invariant under the action
of the structure group of the bundle $\pi : J^k(N,P) \to N \times P$,
the set $\Sigma^I \subset J^k(N,P)$ is well-defined. Clearly
$\Sigma^I$ is an embedded submanifold of $J^k(N,P)$ whose codimension 
is given by \fref{cod}. Given a smooth function $f : N \to P$
whose $k$-jet extension intersects $\Sigma^I$ nontrivially and
transversely at $x_o$,~let   
$$
\Sigma^I_f \stackrel{\dim{def.}}{=} j^kf^{-1} (\Sigma^I)\;.
$$
$\Sigma^I_f$ is a submanifold of $N$ in a \nd of $x_o$, and the
following holds.  

\begin{prop}[\cite{Boardman}, Theorem (6.2)]\label{Sarabande}
$$
{\rm dim\, Ker\,} d(f|_{\Sigma^{(i_1,\ldots, i_k)}_f})(x_o) = \ell \quad
\Longleftrightarrow \quad j^{k+1}f(x_o) \;\; \in \;\;
\Sigma^{(i_1,\ldots, i_k, \ell)}\, . 
$$  
\end{prop}

Let $M$ be an $m$-dimensional manifold endowed with a foliation $\F$
of dimension $n$ and codimension $q$. The inclusion of a leaf
$F$ into $M$ gives rise to a submersion
\begin{equation}\label{Messian}
r_F : \pi^{-1} (F\times \R) \subset J^k(M,\R) \to J^k(F,\R) : j^kf(x)
\mapsto j^k(f|_{F})(x)\;. 
\end{equation}
For every sequence $I = (n = i_1, \ldots,i_k)$ of integers,~let
$$
\Sigma^I_\F \stackrel{\dim{def.}}{=} \bigcup_{F \in \{{\rm leaves \; of
\;} \F\}} r_F^{-1} \left(\Sigma^I\right) \;.
$$
The set $\Sigma^I_\F$ is an embedded submanifold of $J^k(M,\R)$ whose
codimension coincides with that of $\Sigma^I$ in $J^k(n,1)$. Indeed, 
given a foliated chart $(U,\varphi)$ centered at a point $x$, we~have 
$
\pi^{-1}(U \times \R) \cap \Sigma^I_\F \simeq U \times \R \times 
r^{-1}(\Sigma^I)  
$,
where $r$ is the submersion $J^k(m,1) \to J^k(n,1) : j^kf(0) \mapsto
j^k(f|_{\R^n \times \{0\}})(0)$. \\ 

\begin{dfn}\label{basson}
In $C^\infty (M)$ endowed with the fine $C^\infty$ topology, let $\a$ 
denote the set of functions whose $k$-jet extension is transverse to 
$\Sigma^I_\F$ for all $I = (\sq i k)$ and all $k = 1, 2, \ldots\,$. 
\end{dfn}
By Thom's Transversality Theorem, this set is residual. Now given
$f$ in $\a$ and any sequence $I$ of length $k$, the set
$\Sigma^I_{\F,f} = (j^k f)^{-1} (\Sigma^I_\F)$ is an embedded 
submanifold whose codimension equals that of $\Sigma^I_\F$ in
$J^k(M,\R)$; it is called a {\em Boardman stratum}. There are only
finitely many nonempty Boardman strata, provided we ignore sequences
ending with more than one zero. The last assertion follows from the
observation that if $i_{k+1} > 0$, the codimension of
$\Sigma^{(I,i_{k+1})}$ in $J^{k+1}(n,p)$ is strictly larger than that
of $\Sigma^{I}$ in $J^k(n,p)$ (cf.~\cite{Boardman} p.~47).

The following result is a consequence of \pref{Sarabande}.

\begin{prop}\label{Courante} Let $f$ be an element of $\a$ and let
$x_o$ be a point in $\Sigma^{(i_1, \ldots ,i_k)}_{\F,f}$. Then
\begin{equation}\label{Ravel}
{\rm dim }\,\left(T_{x_o}\F \cap T_{x_o} \Sigma^{(i_1, \ldots
,i_k)}_{\F,f}\right) = \ell \quad \Longleftrightarrow \quad j^{k+1}
f(x_o)\;\; \in \;\; \Sigma^{(i_1, \ldots ,i_k,\ell)}_\F 
\end{equation}
\end{prop}

\Pf Since the problem is local, we may assume that $M$ is
diffeomorphic to a product $F\times V$, with $F=\R^n$, $V=\R^q$, and
where $\F$ corresponds to the ``horizontal'' foliation of $F\times V$
(whose leaves are the $F\times \{v\}$'s). Let $\sq u n$ and $\sq v q$
denote the standard coordinates in $F$ and $V$ respectively. To the
smooth function $f : M \to \R$, we associate a smooth map $\ti f : M
\to \R \times V$ by setting $\ti f(x) = (f(x),p(x))$, where $p$
denotes the natural projection $F\times V \to V$. The proof of
\pref{Courante} relies on the following three statements. Let $x$ be
any point in $M$, and let $k$ be any positive integer, then

\begin{enumerate}
\item[a)] $j^kf(x) \in \Sigma^{I}_\F \subset J^k(M,\R)$ \Iff $j^k\ti
          f(x) \in \Sigma^{I} \subset J^k(M,\R \times V)$;
\item[b)] $j^kf$ intersects $\Sigma^{I}_\F$ transversely at $x$ \Iff
          $j^k\ti f$ intersects $\Sigma^{I}$ transversely at $x$;
\item[c)] if $x$ belongs to $\Sigma^{I}_{\F,f}$, then $T_{x}\F \cap T_{x}
          \Sigma^{I}_{\F,f} = {\rm Ker }\,d (\ti
          f|_{\Sigma^{I}_{\F,f}})(x)$.   
\end{enumerate}
Assuming these statements proved, we conclude the proof of
\pref{Courante} as follows. Statements $a)$ and $b)$, combined 
with \pref{Sarabande} imply that  
\begin{equation}\label{Debussy}
{\rm dim\, Ker\,} d(\ti f|_{\Sigma^I_{\F,f}})(x_o) = \ell
\quad \Longleftrightarrow \quad j^{k+1} f (x_o)
\in\Sigma^{(i_1,\ldots, i_k, \ell)}_\F\, . 
\end{equation}
Statement $c)$ implies that (\ref{Debussy}) is equivalent to~(\ref{Ravel}).
\cqfd

\vspace{.2cm}
{\em Proof of $a)$.} Consider the map 
$$\begin{array}{rcl}
J^k(M,\R \times V) & \stackrel{p_1 \times p_2}{\longrightarrow} & 
J^k(M,\R) \times J^k(M,V) \\  
j^kg(x) & \mapsto & \left(j^kg_1 (x) , j^kg_2 (x)\right)\;,
\end{array}$$
where $g$ is a map $M \to \R \times V$, whose components are $g_1 : M
\to \R$ and $g_2 : M \to V$. Let $\p \subset J^k(M,V)$ be the image of
the jet extension $j^kp$, and let $\ti \p = p_2^{-1}(\p)$. It is
sufficient to prove that in $J^k(M,\R \times V)$,  

\begin{equation}\label{Moussorgsky}
\ti \p \cap \Sigma^I = \ti \p \cap p_1^{-1}
(\Sigma^I_\F) \,. 
\end{equation}
Let $z = (z_1,z_2) = (j^k f(x), j^k p (x)) = j^k\ti f(x)$ be an
element in $\ti \p$. We assume for convenience that $x=0$. Let $\ti
z_1$ denote $j^k(f|_{F \times \{0\}}) (x)$. We will prove, by induction 
on $0 \leq \ell \leq k$, that in $\e_m$, 
 
\begin{equation}\label{Clara}
\delta^\ell (\I(z)) = (v) + \delta^\ell(\I(\ti z_1)) \;,
\end{equation}
where $(v)$ denotes the ideal generated by the germs of the functions
$v_1, \ldots, v_q$, and where $\delta^\ell(\I(\ti z_1))$ is thought of
as being an ideal of $\e_m$. This will imply that ${\rm rk}\,
\delta^\ell(\I(z)) = q + {\rm rk}\, \delta^\ell( \I(\ti z_1))$ for all
$0 \leq \ell \leq k$, hence that $z_1$ lies in $\Sigma^I_\F$ \Iff $z$ lies in 
$\Sigma^I$, which is the content of (\ref{Moussorgsky}). \\  

Observe that (\ref{Clara}) clearly holds when $\ell=0$. Indeed, in a
\nd of $0$ in $M$, the function $f$ can be written as    
$$
f(u,v) = f(u,0) + \sum_{i=1}^q v_ih_i(u,v)
$$ 
for some smooth functions $h_i$. Then suppose that (\ref{Clara}) holds
for some $\ell\geq 0$. The ideal $\delta^{\ell+1} (\I(z))$ is
generated by  $\delta^\ell (\I(z))$ and the set $\Gamma_{r + 1}
\delta^\ell (\I(z))$ of $(r+1)\times (r+1)$ minors of the Jacobian
matrix $(D_i f_j)_{1\leq i\leq m,1\leq j\leq a}$, where $r = {\rm
rk}\, \delta^\ell(\I(z))$, and where $f_1, \ldots, f_a$ is a set of 
generators for $\delta^\ell (\I(z))$. Since we may assume that $f_1 = 
v_1, \ldots, f_q = v_q$ and that $f_{q+1}, \ldots, f_a$ are generators of
$\delta^\ell (\I(\ti z_1))$, the set $\Gamma_{r + 1} \delta^\ell (\I(z))$ is
also the set of $(r-q+1) \times (r-q+1)$ minors of $(D_i f_j)_{1\leq i
\leq n, q+1\leq j\leq a}$. Since the latter set coincides with 
$\Gamma_{r-q+1} \delta^\ell (\I(\ti z_1))$, we have 
$$
\delta^{\ell+1} (\I(z)) = (v) + \delta^\ell (\I(\ti z_1)) +
\Gamma_{r-q+1} \delta^\ell (\I(\ti z_1)) = (v) + \delta^{\ell+1}
(\I(\ti z_1))\;.
$$\cqfd
 
\vspace{.2cm}
{\em Proof of $b)$.} Suppose that $j^kf(x) \in \Sigma^I_\F$. Observe that 
the identity (\ref{Moussorgsky}) implies that $j^kf$ intersects 
$\Sigma^I_\F$ transversely at $x$ \Iff 

$$\begin{array}{ccl}
T_z \ti \p & = & T_z \ti \p \cap T_z \left(p_1^{-1}\Sigma^I_\F\right)
+ (j^k \ti f)_{*_x}(T_x M) \\
& = & T_z \left(\ti \p \cap p_1^{-1} \Sigma^I_\F \right)  + (j^k \ti
f)_{*_x}(T_x M)\\ 
& = & T_z \left(\ti \p \cap \Sigma^I\right) + (j^k \ti f)_{*_x}(T_x M)\;,
\end{array}$$ 
where $z = j^k\ti f (x) = (j^kf(x), j^k p(x)) = (z_1,z_2)$. So, in
order to prove $b)$, we need to show that $\Sigma^I$ and $\ti \p$
intersect transversely at $z$, or equivalently that
$$
T_{z_2} \p + (p_2)_{*_z}(T_z \Sigma^I) = T_{z_2}J^k(M,V)\;.
$$
Since $J^k(M,V) \simeq M \times V \times J^k(m,q)$, and since $\p$ is
transverse to $\{x\} \times V \times J^k(m,q)$, it is sufficient to
prove that $T_{p(x)}V$ and $T_{z_2} J^k(m,q)$ are contained in
$(p_2)_{*_z}(T_z \Sigma^I)$. \\

Observe that if $\varphi^t_X$ denotes the local flow of a vector field
$X$ defined on $M$, vanishing at $x$, then $j^k(\ti f
\circ\varphi^t_X)(x)$ belongs to $\Sigma^I$ for all~$t$. Hence    
$$
\xi_X \stackrel{\dim{def.}}{=} \fr{d \left( j^k(\ti f \circ
\varphi^t_X)(x)\right)}{dt}(0) \quad \in \quad T_z \Sigma^I\;.
$$ 
Moreover, the set of vectors $(p_2)_{*_z} (\xi_X)$, as $X$ varies among 
vector fields vanishing at $x$, coincides with $T_{z_2} J^k(m,q)$. So it 
remains to prove that $T_{p(x)}V$ is contained in $(p_2)_{*_z}(T_z 
\Sigma^I)$ as well. This is easily seen once we notice that 
$T_z\Sigma^I$ contains $T_{\ti f(x)} (\R\times V)$, and hence that 
$(p_2)_{*_z}(T_z \Sigma^I)$ contains $(p_2)_{*_z} (T_{\ti f(x)} (\R\times 
V)) = T_{p(x)} V$.      
\cqfd

\vspace{.2cm}
{\em Proof of $c)$.} From $a)$ we know that $\Sigma^{I}_{\F,f} =
\Sigma^I_{\ti f}$. Besides, if $x$ lies in $\Sigma^{I}_{\F,f}$, then

$$\begin{array}{ccl}
{\rm Ker }\,d (\ti f|_{\Sigma^{I}_{\ti f}})(x) & = &
{\rm Ker}\, d\ti f(x) \; \cap \; T_{x} \Sigma^{I}_{\ti f} \\
& = & T_{x} \F \; \cap \; T_{x
} \Sigma^{I}_{\ti f} \;.
\end{array}$$ 
\cqfd

\vspace{.2cm}
Let $f$ be an element of $\a$. We would like to interpret the first
few $\Sigma^I_{\F,f}$'s. Before doing so, we introduce the {\em foliated
second differential of $f$} ({\wrt the foliation $\F$}), that is, the map 
\begin{equation}\label{secdiff}
\ol d{}^2\!f : T_\Sigma \F \oplus T_\Sigma \F \to \R : (X,Y) \to Y(\ti
X(f))\;, 
\end{equation}
where $\ti X$ is any local section of $T\F$ extending $X$. As in the
nonfoliated case, it is easily checked that $\ol d{}^2\!f$ is
well-defined and symmetric.

\begin{enumerate}
\item[-] $\Sigma^{(n)}_{\F,f} = \Sigma_f$ is the foliated singular locus of
         $f$ (according to \dref{lwcr}).
\item[-] $\Sigma^{(n, 0)}_{\F,f}$ is the set of {\em leafwise nondegenerate 
         critical points}, i.e.~the set of \lw critical points $x$ for
         which $\ol d{}^2\!f(x) : T_x\F \oplus T_x\F \to \R$ is
         nondegenerate.   
\item[-] $\Sigma^{(n,i)}_{\F,f}$ is the set of \lw critical points
         $x$ for which $\ol d{}^2\!f(x)$ has rank~$n-i$. 
\item[-] In general, as implied by \pref{Courante}, $\Sigma^{(\sq i
         k)}_{\F,f}$ is the set of points $x$ in $\Sigma^{(\sq i 
         {k-1})}_{\F,f}$ where $T_x\F \cap T_x \Sigma^{(\sq i 
         {k-1})}_{\F,f} = i_k$.  
\end{enumerate}
Observe that $M$ is partitioned into the embedded submanifolds
$\Sigma^I_{\F,f}$, where $I=(i_1,\ldots,i_k,0)$ and $i_k\neq 0$. They
intersect $\F$ transversely, in the sense that $T_x\Sigma^I_{\F,f}
\cap T_x\F = \{0\}$ for all $x$ in $\Sigma^I_{\F,f}$. Finally, the
submanifold $\Sigma^I_{\F,f}$ is of course not closed in general. In
fact,     
$$
\ol{\Sigma^{(\sq i k)}_{\F,f}} \;\; \subset \;\; \bigcup_{(j_1,\ldots
,j_\ell) \geq (i_1,\ldots ,i_k)} \Sigma^{(j_1,\ldots
,j_\ell)}_{\F,f}\;\;\;, 
$$
where the symbol $\geq$ refers to the lexicographical order on the set
of tuples of nonnegative integers.

\begin{rmk}\label{Mendelshon}{\rm A foliated manifold does not
generally admit functions with no leafwise {\em degenerate} critical
points. Indeed, for codimension one foliations, the foliated 
singular locus of such a function is a (not necessarily connected)
closed transversal intersecting every compact leaf. The Reeb foliation
on $S^3$, for example (cf.~\cite{C-C} p.~93), has no closed
transversals intersecting the torus leaf (cf.~\cite{C-C} p.~147). 
}\end{rmk}

\subsection{Isolatedness of leafwise critical points}\label{Milhaud}

We will need (cf.~\oref{piano}) to approximate an $\F$-generic
function by one whose leafwise critical points are isolated \wrt the
leaf topology. The purpose of this subsection is to prove that the
latter property is generic (cf.~\pref{isolate}).\\    

Recall that $\e_n$ and $\frak{m}_n$ denote the ring of germs of maps
$\R^n \to \R$ at the origin and its maximal ideal respectively. The 
set of $k$-jets of maps $(\R^n,0) \to (\R,0)$ at the origin is 
denoted by $J^k(n)$. For $\ell > k$, the natural map $J^\ell(n) \to 
J^k(n) : j^\ell f(0) \to j^k f(0)$ is denoted by $\pi^\ell_k$. As
before, $\sq x n$ are the standard coordinates on $\R^n$. 

\begin{dfn}The {\em Jacobian ideal} of an element $f\in \frak{m}^2_n$, 
denoted by $\J(f)$, is the ideal of $\e_n$ generated by the partial
derivatives $\fr{\partial f}{\partial x_1}, \ldots , \fr{\partial
f}{\partial x_n}$ of $f$. We say that {\em $f$ has codimension $k$} if 
$\dim{cod} (\J(f)) = {\rm dim}_\R\,(\frak{m}_n/\J(f)) = k$. 
\end{dfn}
We will think of a singular element $z = j^kf(0)$ of $J^k(n)$ 
(i.e.~$df(0) = 0$) as a polynomial function on $\R^n$, and   
define its Jacobian ideal and codimension accordingly. \\

\noindent
The following is a consequence of Nakayama's lemma. A proof can be found 
in \cite{Igusa} (Appendix~A.2., Proposition~2.2.). 

\begin{lem}\label{string}If $f - f(0) \in \frak{m}_n^2$ has finite 
codimension, then $0$ is an isolated singularity of~$f$. 
\end{lem}

\begin{dfn}
Two germs $f$ and $g$ in $\e_n$ are said to be {\em equivalent} if there 
exists a germ $\phi$ of local diffeomorphism $(\R^n,0) \to (\R^n,0)$
of $\R^n$ at the origin for which $f = g \circ \phi$.
\end{dfn}
Observe that equivalent germs  have same codimension.

\begin{dfn}
A germ $f \in \e_n$ is {\em $k$-determined} if for any element $g$ in 
$\frak{m}_n^{k+1}$, the germ $f + g$ is equivalent to $f$.
\end{dfn}

\begin{thm}[\cite{JM-68}]\label{epinette} If the codimension of $f
\in \e_n$ is $k$, then $f$ is $(k+2)$-determined. \typeout{see also 
\cite{JMt} p.~27)}
\end{thm}
Consider in $J^k(n)$ the set $Z^k$ consisting of singular $k$-jets 
having codimension strictly larger than $k-2$.

\begin{prop}\label{cornemuse}
The set $Z^k$ is a real algebraic subset of the Euclidean
space $J^k(n)$. 
\end{prop}

Before proving \pref{cornemuse}, let us mention some facts about real 
algebraic sets. Recall that a {\em real algebraic} subset of a
Euclidean space $\R^p$ is the zero locus of a collection of polynomial
functions defined on $\R^p$. Such a set is a {\em variety} when it is
irreducible, that is, when it cannot be decomposed into the union of
two proper real algebraic sets. A real algebraic subset $S$ of $\R^p$
can be written uniquely (up to reordering the factors) as a finite
union      
$$
S = V_1 \cup \ldots \cup V_n
$$
of varieties, where no factor $V_i$ is contained in a factor $V_j$ 
with $i \neq j$. Such a decomposition is said to be {\em minimal}. Real 
algebraic sets can also be decomposed into smooth submanifolds~:    

\begin{thm}[\cite{Whitney}]\label{Whitney} Any real algebraic subset
of $\R^p$ can be decomposed into a finite union of disjoint embedded
submanifolds of~$\R^p$.  
\end{thm}
The {\em dimension of $S$} is defined to be the maximum dimension of 
the factors. The following result will be crucial later on (cf.~proof of 
\pref{Schnittke}).   

\begin{thm}[\cite{Whitney}]\label{Janacek}
A proper subvariety of a variety has strictly smaller dimension.
\end{thm} 

{\em Proof of \pref{cornemuse}.} This type of argument is quite standard
(see \cite{Igusa}). Let $f \in J^k(n)$ be a singular jet. First observe 
that $\dim{cod} (\J(f)) > k-2$ \Iff $\dim{cod}(\J(f) + \frak{m}_n^k) > 
k-2$. Indeed, suppose that $\dim{cod} (\J(f) + \frak{m}_n^k) \leq k-2$. 
Then, in the following nested sequence of spaces
$$
\frak{m}_n \supset \J(f) + \frak{m}_n^2 \supset \J(f) + \frak{m}_n ^3 
\supset \ldots  \supset  \J(f) + \frak{m}_n^k\;, 
$$
equality must occur somewhere, that is, $\J(f) + \frak{m}_n ^{r} = 
\J(f) + \frak{m}_n^{r+1}$ for some $r\leq k-1$. This implies that 
$\frak{m}_n^r \subset \J(f) + \frak{m}_n ^{r+1}$. Hence
$\frak{m}_n^{k} \subset \frak{m}_n^{r} \subset \J(f)$ (cf.~\cite{JMt}
Proposition~4.~p.~3). Therefore $\dim{cod}(\J(f)) = \dim{cod}(\J(f) +
\frak{m}_n^k) \leq k-2$. The other implication is obvious.\\

Now $\dim{cod}(\J(f) + \frak{m}_n^k) > k-2$ is equivalent to 
$$
\dim{dim} (\J(f) + \frak{m}_n^k )/\frak{m}_n^k < \dim{dim}(\frak{m}_n/
\frak{m}_n^k) - (k-2) \stackrel{\dim{def.}}{=} N\;.
$$
The vector space $(\J(f) + \frak{m}_n^k )/\frak{m}_n^k$ is generated by 
elements of the type $x^\alpha \frac{\partial f}{\partial x^i}$, where 
$\alpha$ is a multi-index with $|\alpha| \leq k - 2$. Therefore, its 
dimension is less than $N$ when the rank of a certain matrix $M$, whose 
coefficients are linear functions of the coefficients of the
polynomial function $f$, is bounded by $N-1$; or equivalently, when
all minors of $M$ of size $N$ and larger vanish.
\cqfd

\vspace{.2cm}
The combination of \tref{cornemuse} and \tref{Whitney} implies that
$$
Z^k = M_1 \coprod \ldots \coprod M_r\;,
$$
where each $M_i$ is an embedded submanifold of $J^k(n)$. 

\begin{prop}
For all $\ell > k$, the real algebraic set $Z^\ell$ is a subset of the real 
algebraic set $(\pi^\ell_k)^{-1}(Z^k)$.
\end{prop}

\Pf Let $f$ be an $\ell$-jet whose truncation $\pi^{\ell}_k(f)$ does not 
belong to $Z^k$, that is, $\dim{cod} (\J(\pi^{\ell}_k(f))) \leq k-2$. 
\tref{epinette} implies that $f$ is equivalent to $\pi^{\ell}_k(f)$. Hence 
$\dim{cod} (\J(f)) \leq k-2$ as well. Thus $f$ lies outside $Z^\ell$.
\cqfd

\begin{prop}\label{Schnittke}
The codimension $c_k$ of $Z^k$ in $J^k(n)$ is an unbounded function of 
$k$. 
\end{prop}

\Pf Consider the homogeneous polynomial $r = x_1^{k+1} + \ldots + 
x_n^{k+1}$ of degree $k+1$. The Jacobian ideal of $r$ coincides with 
$(x_1^k, \ldots, x_n^k)$, and contains $\frak{m}_n^{n(k-1)+1}$. Thus 
$$
\dim{cod} (r) \leq \dim{cod} (\frak{m}_n^{n(k-1)+1})
\stackrel{\dim{def.}}{=} a_k \;.
$$
Let $p$ be any polynomial of degree $k$. The codimension of $p+r$ is 
bounded by the codimension of $r$ (cf.~\cite{Igusa} A.2., Theorem~2.7).
Thus $\dim{cod} (p+r) \leq a_k$ for all $p \in J^k(n)$. Let $\ell = 
a_k + 2$. We have just shown that for all $p \in Z^k$, the set 
$(\pi_k^{\ell})^{-1} (p)$ is not entirely contained in $Z^{\ell}$.\\

Now, let $Z^k = V_1^k \cup \ldots \cup V^k_{n_k}$ be the minimal 
decomposition of $Z^k$ into varieties. Let $(\pi_k^{\ell})^{-1}(V_j^k)$ 
be denoted by $W_j$. The real algebraic set $W_j$ is irreducible as
well. Because $Z^{\ell} = V_1^\ell \cup \ldots \cup V^\ell_{n_\ell}
\subset (\pi_k^{\ell})^{-1}(Z^k)$, the set $V^{\ell}_i$ coincides with
$(V^{\ell}_i \cap W_1) \cup \ldots \cup (V^{\ell}_i \cap
W_{n_k})$. Since $V^{\ell}_i$ is irreducible, there exists an index
$j_i$ such that $V^{\ell}_i = V^{\ell}_i \cap W_{j_i}$. The previous
paragraph implies that $V^{\ell}_i$ is a proper subset of
$W_{j_i}$. Hence, the codimension of $V^{\ell}_i$ is strictly larger
than that of $W_{j_i}$ (\tref{Janacek}). Thus, the codimension of
$Z^{\ell}$ in $J^{\ell}(n)$ is strictly larger than that of
$(\pi^{\ell}_k)^{-1} (Z^k)$, that is $\dim{cod} Z^{\ell} > \dim{cod}
Z^k$.  
\cqfd 

\begin{prop}\label{isolate}Let $M$ be a manifold endowed with a
foliation $\F$. For $f$ in a residual subset $\b$ of $C^\infty(M)$, and for
any leaf $F$ of $\F$, the critical points of $f|_F$ are isolated in~$F$.
\end{prop}

\Pf Given a leaf $F$ of the foliation $\F$, define $Z^k(F) \subset 
J^k(F,\R)$ to be   
$$
Z^k(F) = \bigcup_{(x,y)\in F \times \R} Z^k(F)_{(x,y)}\;,
$$ 
where $Z^k(F)_{(x,y)}$ corresponds to $Z^k$ under the identification of 
the fiber at $(x,y)$ of the bundle $J^k(F,\R) \to F \times \R$ with
$J^k(n)$. Then, define $Z^k(\F) \subset J^k(M,\R)$ to~be~:   
$$
Z^k(\F) = \bigcup_{F \in \{ {\rm leaves \; of \;} \F \} } r_F^{-1}(Z^k(F))\;,
$$
where $r_F$ is defined by the expression (\ref{Messian}). Locally, the set 
$Z^k(\F)$ is the union of finitely many embedded submanifolds of $J^k(M,\R)$ 
having codimension at least $c_k$. Hence, Thom's Transversality Theorem 
implies that, provided $k$ is large enough for $c_k > {\rm dim}\, M$, the 
set $\b$ of smooth real-valued functions on $M$ whose $k$-jet does not meet 
$Z^k(\F)$ is residual in $C^\infty(M)$ \wrt the fine $C^\infty$ topology. 
Let $f$ be such a function. For any \lw singularity $x$ of $f$, the jet 
$j^k(f|_{F_x})(x)$ has codimension $k-2$ at most. Thus, the germ of 
$f|_{F_x}$ at $x$ is equivalent to $j^k(f|_{F_x}) (x)$, and has finite 
codimension as well. \lref{string} implies that the singularities of 
$f|_{F_x}$ are isolated in~$F_x$.          
\cqfd

\subsection{$\F$-generic functions}

Let $M$ be an $m$-dimensional manifold carrying a foliation $\F$ of
dimension $n$ and codimension~$q$. 

\begin{dfn}\label{cello}($\F$-genericity.) A smooth real-valued
function $f$ on $M$ is said to be {\em $\F$-generic} if $j^1f$
intersects $\Sigma^{(n)}_{\F}$ transversely and $j^2f$ intersects
$\Sigma^{(i_1,i_2)}_{\F}$ transversely for all pairs~$(i_1,i_2)$.  
\end{dfn}

\begin{dfn}\label{strongcello}(Strong $\F$-genericity). A smooth
real-valued function on $M$ is said to be {\em strongly
$\F$-generic}~if    
\begin{enumerate}
\item[i)] the critical points of $f$ are nondegenerate and
          \lw nondegenerate;
\item[ii)] for every sequence $I$ of any length $k$, the map $j^kf$
           intersects $\Sigma^I_{\F}$ transversely;  
\item[iii)] the \lw critical points of $f$ are isolated \wrt the leaf
            topology.
\end{enumerate}
\end{dfn}

Here as well, it is a consequence of Thom's Transversality Theorem that
the set of (strongly) $\F$-generic functions is residual in
$C^{\infty}(M)$ for the fine $C^{\infty}$ topology. Indeed, the set of
functions whose singularities are nondegenerate is a residual set
(cf.~\cite{Milnor}), as is the set of functions satisfying $ii)$ 
(introduced in \dref{basson}), and the set of functions satisfying $iii)$
(containing the set $\b$ introduced in \pref{isolate}). So the only
thing that remains to be proven is that generically, critical points of 
$f$ are \lw nondegenerate. Introduce the following subsets of~$J^2(M,\R)$
$$
W^{\ell} \stackrel{\dim{def.}}{=} \Sigma^{(m,0)} \cap
\Sigma^{(n,\ell)}_\F\quad 1 \leq \ell \leq n \;. 
$$
One easily proves that the $W^\ell$'s are embedded submanifolds of
$J^2(M,\R)$ of codimension strictly larger than $m$. Hence, the set of
functions whose second jet does not meet any $W^\ell$ is residual. By
construction, the critical points of those functions are \lw
nondegenerate. 

\begin{rmk}\label{Verdi}{\rm Reflecting on \dref{Mahler}, one might
think that $c)$ is superfluous. Indeed, the set of strongly
$\F$-generic functions is dense in $C^{\infty}(M)$, and we can
therefore approximate in the fine $C^{\infty}$ topology any given
function by a strongly $\F$-generic one. The problem is that one
usually creates leafwise local maxima in the process. It is not clear,
although quite plausible, that a function $f$ with no leafwise local
maxima always admits a nearby strongly $\F$-generic function with no
\lw local maxima. However, if $f$ is already $\F$-generic, then any
sufficiently ($C^{\infty}$-) nearby strongly $\F$-generic function will
not have leafwise local maxima either.}   
\end{rmk} 

The proof of the next lemma is very similar to that of the corresponding 
result in the single-leaf case (cf.~\cite{Milnor} for example). Nevertheless, 
we include the proof for the sake of completeness.    

\begin{lem}\label{tuba}
Let $f : M \to \R$ be a smooth function, and let $x$ be a point in $M$
for which $j^2f(x) \in \Sigma{}^{(n,n-r)}_\F$. Then there exists a
chart centered at $x$ and adapted to $\F$, with local coordinates
$x_{1},\ldots,x_{m}$, \wrt which $f$ has the following expression  
\begin{equation}\label{virginal}
f(x_{1},\ldots x_{m}) = \pm x_{1}^{2} \pm \ldots \pm x_{r}^{2} + 
f(0,\ldots,0,x_{r+1},\ldots ,x_{m})\;,
\end{equation}
\end{lem}

\Pf
Consider local coordinates $x_{1},\ldots,x_{m}$ centered at $x$ for
which $\F$ is defined near $x$ by the equations $x_{n+1}=c_{n+1},
\ldots,x_{m}=c_{m}$, where $c_1,\ldots,c_m$ are constants . Since 
$j^2f(x) \in \Sigma{}^{(n,n-r)}_\F$, the function $f|_{F_x}$ is 
singular at $x$, and we may assume (after eventually performing a 
linear change of the coordinates $\sq x n$) that 

$$\begin{array}{ccc}
\displaystyle{\left(\frac{\partial^{2} f}{\partial x_{i}\partial 
x_{j}}(0)\right)_{1\leq i,j\leq r}} & = & \left(\begin{array}{cc}
                                      I_d & 0 \\
                                      0     & -I_{r-d}
                                      \end{array}\right)\;.

\end{array}$$
One deduces from the implicit function theorem that there exists a \nd
$U$ of $0$ in $\R^{m-r}$, a \nd $U'$ of $0$ in $\R^{r}$, and a function 
$\alpha : (U,0)\to (U',0)$ such that, for every $i =1,\ldots, r$, and for
$x$ in $U' \times U$, one~has 

$$\begin{array}{ccc}
\displaystyle{\frac{\partial f}{\partial x_{i}}(x_{1},\ldots,x_{m})} = 0 
& \mbox{\Iff} & (x_{1},\ldots,x_{r}) = \alpha(x_{r+1}, \ldots ,x_m)\;.
\end{array}$$
Use $\alpha$ to define new coordinates near $x$ as follows~:
$$\left\{\begin{array}{cclll}
x'_{i} & = & x_{i} - \alpha_{i}(x_{r+1},\ldots,x_{m}) & 
\mbox{for} & 1\leq i\leq r \\
x'_{j} & = & x_{j} & \mbox{for} & j>r\;.
\end{array}\right.$$
Now 

$$
\frac{\partial f}{\partial 
x'_{i}}(0,\ldots,0,x'_{r+1},\ldots,x'_{m}) = 0\quad \mbox{for} 
\quad 1\leq i\leq r \;.
$$ 
The primes are omitted hereafter. A classical argument
(cf.~\cite{Milnor} Lemma 2.1.) shows that in a convex \nd of~$0$,
one~has   

\begin{equation}\label{4}
f(x_{1},\ldots,x_{m}) = f(0,\ldots,0,x_{r+1},\ldots,x_{m}) + 
\sum_{1\leq i\leq j\leq r} x_{i}x_{j}f_{ij}(x_{1},\ldots,x_{m})\;.
\end{equation}
Consider the coordinates $x'_1,\ldots, x'_m$ defined near $x$~by 

$$\left\{\begin{array}{ccl}
x'_{1} & = & \displaystyle{x_{1}\sqrt{|f_{11}(x)|} + \fr 1 2
\sum_{i=2}^{r} x_{i}\frac{f_{1i}(x)}{\sqrt{|f_{11}(x)|}}} \\
x'_{i} & = & x_{i}\quad\mbox{for}\quad i\geq 2\;.
\end {array}\right.$$
With respect to these coordinates, $f$ can be written as (\ref{4}),
where $f_{11}(x) = 1$, where $f_{1i}(x) = 0$, and where $f_{ij}(0)$
remains unchanged. We repeat this process, modifying successively
$x_2, \ldots ,x_r$, until the expression (\ref{virginal}) is
achieved. Observe that every change of coordinates that has been made
preserves $\F$ since it leaves the coordinates $x_{n+1},\ldots,x_m$
untouched.  
\cqfd

\begin{prop}\label{positivity}
Let $f:M \to \R$ be an $\F$-generic function without leafwise local 
maxima, and let $x$ be in $\Sigma^{(n,n-r)}_{\F,f}$. Then, for any
system of coordinates centered at $x$ satisfying the conclusion of the
previous lemma, the quadratic form $Q(x_{1},\ldots,x_{r})= \pm
x_{1}^{2} \pm \ldots\pm x_{r}^{2}$ occurring in expression (\ref{virginal}) 
must have at least one positive~sign.  
\end{prop}

\Pf
Suppose on the contrary that either $Q=0$ or $Q$ is negative definite. 
We claim that, if so, any neighborhood of $x$ contains leafwise local 
maxima of $f$, contradicting the assumption.\\

Let $x\in M$, and let $(U,\varphi)$ be a chart at $x$, adapted to 
$\F$, with local coordinates $x_{1}, \ldots,x_{m}$, \wrt which $f$ is 
expressed by (\ref{virginal}). The chart $(U,\varphi)$ induces a chart
$(\ti{U},\ti{\varphi})$ on $J^{2}(M,\R)$~:

$$\begin{array}{ccccl}
\ti{\varphi} & : & \ti{U} = (\pi^{2})^{-1}(U) & \to &\varphi(U)\times
\R\times\R^{m}\times {\mathbb S}^m \\ 
&& p = j^{2}g(x) & \mapsto & (\varphi(x), g(x), L_p, S_p)\;,
\end{array}$$
where $\pi^2$ denotes the natural projection $J^2(M, \R) \to M$, and
where ${\mathbb S}^m$ denotes the set of symmetric $m\times m$ matrices 
with real coefficients. The symbol $L_p$ (\rp $S_p$) stands for the 
Jacobian (\rp Hessian) matrix of $g \circ \varphi^{-1}$ at
$\varphi(x)$. \\

It is convenient to decompose a symmetric $n\times n$ matrix $S$ into
blocks as follows~:     
$$
S = \left(\begin{array}{cc} 
                         A & B \\
                         B^{t} & D
                         \end{array}\right)\;,
$$
where $A$ and $D$ are $r\times r$ and $(n-r)\times(n-r)$ symmetric 
matrices, and where $B$ is some $r\times(n-r)$ matrix. If $A$ is 
nonsingular, then $S$ has rank $r$ \Iff the matrix $e(S) =
-B^{t}A^{-1}B + D$ vanishes. Define ${\cal U}$ to be the neighborhood
of $j^2f(x)$ in $\ti{U}$ consisting of jets of local maps whose
Hessian has a nonsingular upper left $r\times r$ block. Thus 
$$
\Sigma^{(n,n-r)}_{\F,f}\cap{\cal U} = (\rho \circ \ti{\varphi} 
)^{-1}(0)\;,
$$ 
where $\rho$ denotes the submersion~:    
$$\begin{array}{cclll}
\rho = \rho^1 \times \rho^2 & : & \ti{\varphi}({\cal U}) & \to &
\R^{n}\times {\mathbb S}^{n-r} \\ 
{} & {} & (x,t,L,S) & \mapsto & \left( p(L), e(q(S))\right)\;\;,
\end{array}$$

\noindent and where $p : \R^n\times \R^{m-n} \to \R^n$ (\rp $q :
{\mathbb S}^m \to {\mathbb S}^n$) is the natural projection. Since
$j^{2}f$ (respectively $j^{1}f$) intersects $\Sigma^{(n,n-r)}_\F$
(respectively $\Sigma_f$) transversely at $x$, the maps $\rho \circ
\ti{\varphi} \circ j^{2}f$ and $\rho^1 \circ \ti{\varphi} \circ
j^{2}f$ are submersions near $x$. Moreover $(\rho^1 \circ \ti{\varphi}
\circ j^{2}f)_{*_x} (T_x \Sigma_f) = 0$, 
and~thus

$$
\rho^{2}\circ\ti{\varphi}\circ j^{2}f|_{\Sigma_{f}\cap U} : 
\Sigma_{f}\cap U  \to  {\mathbb S}^{n-r} 
$$
is a submersion near $x$ as well. This implies in particular that for 
$\varepsilon$ sufficiently small, there exists a $x_o$ in $\Sigma_{f}
\cap U$ such that $\rho^{2} \circ \ti{\varphi} \circ j^{2}f (x_o) = 
-\varepsilon I_{n-r}$. Then $\ol{d}{}^2\!f (x_o) = - \dim{Id}$ is
negative definite. Thus $f$ achieves a leafwise local maximum at $x_o$.
\cqfd

\begin{rmk}{\rm The assumption that $f$ is $\F$-generic and has no
leafwise local maxima prevents $f$ from having totally degenerate
singularities (singularities at which the foliated second differential
of $f$ vanishes).} 
\end{rmk}

\section{Leafwise gradient vector fields}\label{lwgr}

Let $(M,\F)$ be a foliated manifold, let $f$ be a positive, proper, strongly
$\F$-generic function on $(M,\F)$ (\dref{strongcello}), and let $g$ be
a Riemannian metric on $M$. 

\begin{dfn}\label{lwgrdef}The {\em leafwise gradient vector field of
$f$ associated to the metric $g$}, denoted by $\nabla_\F f$, is the
orthogonal projection onto $T\F$ of the (ordinary) gradient vector
field of $f$. Equivalently, $\nabla_\F f$ is the section of $T\F$
whose restriction to each leaf $F$ is the gradient of $f|_{F}$
associated to the metric~$g|_{F}$.
\end{dfn}
The local flow of $\nabla_\F f$ is denoted by $\varphi^{t}$. The orbit of 
a point $x$ under the flow $\varphi^t$ is denoted by $\theta_x$. If 
$S$ is some subset of $\R$, the portion $\{\varphi^t(x) ; t\in S\}$ of 
$\theta_x$ is denoted by $\varphi^S(x)$. Finally, the negative (\rp
positive) limit set  of the orbit $\theta_x$ is denoted by $\EK^-_{x}$
(\rp $\EK^+_{x}$), i.e.
$$
\displaystyle{\EK^-_{x} = \bigcap_{n=0}^{\infty} \ol{\varphi^{(-\infty,
-n)} (x)}} \quad \mbox{and} \quad \displaystyle{\EK^+_{x} =
\bigcap_{n=0}^{\infty} \ol{\varphi^{(n, \infty)} (x)}}\;.
$$

\pagebreak

\begin{lem}\label{Britten}
\hspace{-.5cm}\begin{enumerate}
\item[i)] The set of points where $\nabla_\F f$ vanishes coincides
          with $\Sigma_{f}$. Where $\nabla_\F f$ does not vanish,
          $(\nabla_\F f) (f)>0$. 
\item[ii)] For every $x$ in $M$, the set $\EK^-_{x}$ reduces to a 
single point, lying in $\Sigma_{f}$. Similarly, the set $\EK^+_{x}$
is either empty if the trajectory of $x$ is unbounded, or consists 
also of a single point, lying in $\Sigma_{f}$ as~well. 
\end{enumerate}
\end{lem}

\Pf
\begin{enumerate}
\item[$i)$] Both assertions follow from the observation that
for any $X$ in $T\F$, we~have
$$
Xf = df(X) = g(\nabla f, X) = g(\nabla_\F f, X)\;.
$$ 
 
\item[$ii)$] Since $\varphi^{(-\infty , 0]}(x)$ is bounded 
($\varphi^{(-\infty , 0]}(x) \subset f^{-1} [0,f(x)]$), for any sequence 
$t_{k}\to -\infty$ of real numbers, the sequence $\varphi^{t_{k}}(x)$
has a converging  subsequence . The set $\EK^-_x$ is therefore  
nonempty. Moreover, because $\EK^-_x$ is $\varphi^t$-invariant and 
contained in a level set of $f$ (the one at level ${\rm inf}\,
\{f(\varphi^{t}(x)); t\in \R\}$), it must be contained in $\Sigma_{f}$ 
(since $f$ increases along nonconstant trajectories of~$\nabla_\F f$). \\    

Observe now that $\EK^-_x$ is connected. Indeed, suppose it is not, and
let $O_{1}$ and $O_{2}$ be disjoint open subsets of $M$ containing
$\EK^-_x$ in their union. There exists a $T\leq 0$ for which $\varphi^{t}
(x)\in O_{1}\cup O_{2}$ as soon as $t\leq T$. Otherwise there would be
a sequence $t_{k} \to -\infty$ for which $O_{1} \cup O_{2}$ does not
contain $\{\varphi^{t_{k}}(x)\}$, hence not either $\lim_{k \to \infty}
\varphi^{ t_{k}}(x)$. Connectedness of $\varphi^{(-\infty,T]}(x)$
implies that $\EK^-_{x}$ is entirely contained in one of the~$O_i$'s. \\

To conclude, we need to prove that $\EK^-_x$ is
contained in the leaf $F_x$. Suppose not, and let $w$ be an point in
$\EK^-_x$ that does not belong to $F_x$. Let also $(U,\varphi)$ be a
chart adapted to $\F$ centered at $w$ with $\varphi (U) = [0,1]^n
\times [0,1]^q$. Denote by $P_a$ the plaque $\varphi^{-1}([0,1]^n 
\times \{a\})$. Since $f$ is strongly $\F$-generic, its \lw critical
points are isolated for the leaf topology. We may therefore assume
that $w$ is the only point in $\Sigma_f$ that lies in $P_0$. Because $w$
belongs to $\EK^-_x - F_x$, there is a sequence $a_k$ converging to $0$
such that each plaque $P_{a_k}$ contains a segment $\varphi^{[t_k,
s_k]}(x)$ of $\theta_x$ that meets the boundary of $P_{a_k}$ along its
endpoints (i.e.~$\varphi^{t_k}(x)$ and $\varphi^{s_k}(x)$ are in
$\partial P_{a_k}$). We may also assume, without loss of generality,
that $t_k$ and $s_k$ both converge to $-\infty$. Now the sequence
$\varphi^{t_k}(x)$ has a subsequence that converges to a point $w'$ in
$\partial P_0$. Since by construction $w'$ lies in $\EK^-_x$, it is
necessarily a \lw critical point of $f$, contradicting the hypothesis
made earlier on~$U$. \\   

We have now all the ingredients needed~: $\EK^-_{x}$ is nonempty and
connected, consists of \lw critical points of $f$, and is entirely
contained in a leaf. Property $iii)$ of \dref{strongcello} implies
thus that $\EK^-_{x}$ is a single point. The same argument proves the
second part of statement $ii)$ as~well.   
\cqfd
\end{enumerate}

\begin{dfn}\label{violoncelle}
Given $x$ in $M$, the limit point $\EK^-_x$ is denoted by $x^-$. 
Similarly, when $\theta_x$ in bounded, the limit point $\EK^+_x$ is 
denoted by $x^+$.
\end{dfn}

\begin{dfn}\label{stable}
\hspace{-.5cm}\begin{enumerate}
\item If $A$ is a subset of $\Sigma_{f}$, we define the {\em stable
      set of $A$} to be the set $\W(A) = \{x \in M\, ;\, x^{+}\in
      A\}$. 
\item If $A$ is a subset of $M$, the {\em saturation} of $A$ is
      defined to be the set $\ES A = \{x\in M\,;\,\varphi^t(x)\in A
      \mbox{ for some } t\geq 0\}$.
\end{enumerate}
\end{dfn}

The following lemma will reduce the proof that the h-principle holds
on $M$ to the proof that it holds near the ``skeleton''
$\W(\Sigma_{f})$. Let $(a, b)$ be any pair of real numbers with
$a<b$. The {\em slice} $f^{-1}([a,b])$ is denoted hereafter by
$M^b_a$, and the set $\W(\Sigma_f \cap M^b_a)\cap M^b_a$ by
$\W^b_a(\Sigma_f)$.  

\begin{lem}\label{alto}
Given any pair $a<b$ of regular values of $f$, and given any
neighborhood $U$ of  $M^a_0 \cup \W^b_a(\Sigma_{f})$, the flow of
$\nabla_\F f$ yields a foliated isotopy  $\psi_{t}: M\to M$, $t \in
[0,1]$, such~that  
\begin{enumerate}
\item[-] $\psi_{t}(M^b_0)\subset M^b_0$ for all $t$ in $[0,1]$,
\item[-] $\psi_{1}(M^b_0)\subset U$,
\item[-] $\psi_{t}$ coincides with the identity map on $M^a_0 \cup 
         \W^b_a(\Sigma_{f})$,
\item[-] $\psi_{t}$ is stationary for $t\geq\frac{1}{2}$.
\end{enumerate}
\end{lem}

\Pf 
Let $U$ be an open neighborhood of $M^a_0 \cup \W^b_a(\Sigma_{f})$.
As a first step, we prove existence of a neighborhood $N$ of
$\Sigma_{f} \cap M^b_a$ whose saturation $\ES N$ is contained in
$U$. Suppose that such an $N$ does not exist. Consider the
neighborhoods $N_{k}$ of $\Sigma_{f} \cap M^b_a$ given by
$N_{k} = \{x\in M ; $ the distance between $x$ and $\Sigma_{f} \cap
M^b_a$ is less than $\frac{1}{k}\}$. Since $\Sigma_{f} \cap M^b_a$ is
compact, $N_{k}$ is in $U$ for $k$ large enough. Now, suppose that for
every $k$ there exists $x_{k}$ in $N_{k}$ such that for some negative
$t_{k}$, the point $y_{k}=\varphi^{t_{k}}(x_{k})$ does not belong to
$U$. There exists subsequences $x_{k_{l}}$ and $y_{k_{l}}$ of $x_{k}$ 
and $y_{k}$ respectively, converging to points $x$ and $y$ respectively. 
The point $x$ lies in $\Sigma_{f} \cap M^b_a$, and the point $y$ lies 
outside $U$. Thus the trajectory $\theta_{y}$ intersects $f^{-1}(b) - 
\Sigma_{f}$ nontrivially. Denote by $\h$ the foliation of $M-\Sigma_{f}$ 
by the flow lines of $\nabla_\F f$. Let $\KK$ be the portion of 
$\theta_{y}$ located between levels $f(y)$ and level $b$. Since $\KK$ 
is a compact subset of a leaf of $\h$, there exists a closed 
neighborhood $O$ of $\KK$ in $M-\Sigma_{f}$ such that $\h|_{O}$ is 
isomorphic to a product foliation. For $l$ sufficiently large, 
$y_{k_{l}}$ is contained in $O$, and therefore $x_{k_{l}}$ and $x$ are 
contained in $O$ as well. This contradicts the hypothesis that $O$ does 
not intersect~$\Sigma_{f}$. \\

Now choose an open neighborhood $N$ of $\Sigma_{f} \cap M^b_a$ for
which $\ES N\subset U$. Notice that $\ES N$ is an open neighborhood of 
$\W^b_a(\Sigma_f)$. Since $f^{-1}(a)$ is compact, there exists a 
$\delta>0$ such that $f^{-1}([0,a+\delta))\subset U$. Let $\tilde U = 
\ES N\cup f^{-1}([0,a+\delta))$. Because $f$ increases along
nonconstant trajectories of $\nabla_\F f$, there exists, for every $x$
in $M^b_a$, a $t_{o}\leq 0$ such that $\varphi^{t}(x) \in \tilde U$
for all $t\leq t_{o}$. Moreover, since $M^b_a$ is compact, such a
$t_{o}$ can be chosen that does not depend on $x$. Hence $\varphi^{t_{o}}$
shrinks $M^b_0$ inside $\tilde U$, but does not coincide with the
identity on $M^a_0 \cup \W^b_a(\Sigma_{f})$. Multiplying $\nabla_\F f$
by a smooth function $\eta : M \to [0,1]$ that equals $0$ on 
$M^a_0 \cup \W^b_a(\Sigma_{f})$ and $1$ outside $\tilde U$ yields a
vector field $\eta \nabla_\F f$ whose flow provides us with the
desired isotopy. Indeed, let $\psi_t = \varphi^{\theta(t)\, t_o}$,
where $\theta$ is a smooth function $[0,1] \to [0,1]$ such that 
$\theta(0) = 0$ and $\theta(t) = 1$ for all $t\geq \fr 1 2$. The
isotopy $\psi_t$ satisfies the four required properties.
\cqfd

\section{Construction of the metric}\label{metric}

Let $M$ be an $m$-dimensional manifold endowed with a foliation $\F$ of
dimension $n$ and codimension $q$. Supposing the \fm $(M, \F)$ open
(\dref{Mahler}), let $f : M \to \R$ be a strongly $\F$-generic
function (\dref{strongcello}) without \lw local maxima
(cf.~\rref{Verdi}). In this section, we carry out the construction of a
Riemannian metric for which the stable set of the foliated singular
locus of $f$ (\dref{lwcr} and \dref{stable}) is locally contained in a
finite union of embedded submanifolds, of codimension at least one,
transverse to the foliation $\F$.\\ 

A submanifold $S$ of $M$ is said to intersect $\F$ transversely if for
all $s$ in $S$, the quantity ${\rm dim }\,(T_sS + T_s\F)$ is maximal,
that is, coincides with ${\rm min} \{ {\rm dim }\,S + {\rm dim }\,\F
,\, {\rm dim}\, M\}$. Observe that a submanifold $S$ whose dimension 
is strictly less that $q$ does never intersect each leaf of $\F$ 
transversely. Nevertheless, it might intersect the foliation $\F$ 
transversely. Let $g$ be a Riemannian metric on $M$. Recall the foliated 
second differential $\ol d{}^2\!f$ of $f$ (\sbref{TBs}, \fref{secdiff}). 
Because the bilinear map $\ol d{}^2\!f(\sigma)$ is symmetric for all 
$\sigma \in \Sigma_f$,
the identity $\ol d{}^2\!f(\sigma)(X,Y) = g(X,AY)$ defines a
$g$-self-adjoint map $A : T_\sigma\F \to T_\sigma\F$. The map $A$ can
be diagonalized by means of a $g$-orthogonal basis. We therefore have
a splitting $T_\sigma\F = V^+_\sigma \oplus V^-_\sigma \oplus
V^0_\sigma$ into positive, negative and null eigenspaces for $A$. The
space $V^0_\sigma$ is the kernel of $\ol d{}^2\!f(\sigma)$ and is thus
intrinsically defined (independently of $g$), unlike the other two. The
distributions $V^{\pm} : \sigma \in \Sigma_f \mapsto V^{\pm}_\sigma$,
although not of constant rank, are {\em smooth} (unlike the distribution
$V^0 : \sigma \in \Sigma_f \mapsto V^0_\sigma$). Smoothness for such a
distribution is defined as follows~: for every $\sigma$ in $\Sigma_f$
there exist local sections $X_1, \ldots, X_d$ of $TM$ defined in a
\nd $U$ of $\sigma$ in $\Sigma_f$ such~that 
\begin{enumerate} 
\item[-] $\{X_1(\sigma), \ldots, X_d(\sigma)\}$ is a basis of 
$V^{\pm}_\sigma$, 
\item[-] $X_i(\sigma')$ belongs to $V^{\pm}_{\sigma'}$ for all
$\sigma'$ in $U$. 
\end{enumerate}
 
\begin{rmk}\label{rmkstratum}{\rm Consider a sequence $I=(\sq{i}{k-1}, 
i_k=0)$. The rank of $V^0$ and hence of $V^+\oplus V^-$ is constant 
on the submanifold $\Sigma^I_{\F,f}$ (cf.~\sref{TBs}). Moreover, since 
both $V^+$ and $V^-$ are smooth distributions, their rank must be 
constant on any connected component of $\Sigma^I_{\F,f}$. 
}\end{rmk}

\begin{dfn}\label{defstratum}
The union of the connected components of $\Sigma^I_{\F,f}$ on which
${\rm rk}\,V^+$ achieves a fixed value $d$ is denoted by
$\Sigma_d^{I}$, and called hereafter a {\em stratum of $\Sigma_f$}
(provided it is nonempty). We introduce an order on the set of strata
as follows~:   
$$
\Sigma_{e}^{(\sq j \ell)} \; \leq \;\Sigma_d^{(i_1,\ldots,i_k)}
\Longleftrightarrow  (n-e, \sq j \ell) \geq  (n-d, i_1,\ldots,i_k)\;,
$$
where $\geq$ refers to the lexicographical order on tuples of nonnegative
integers. 
\end{dfn}

Observe that because $i_k=0$, a stratum $S$ is an embedded submanifold
(generally neither connected nor closed) that intersects $\F$
transversely, in the sense that $T_sS \cap T_s\F=\{0\}$ for all $s$ in
$S$. Moreover, there are only finitely many strata, and they
partition~$\Sigma_f$. Furthermore, the closure of $\Sigma_d^{I}$ is
contained in the union of the $\Sigma_e^{J}$'s with $e\leq d$ and $J
\geq I$ (see end of \sref{TBs}). In particular, it is
contained in the union of the smaller (\wrt the order introduced in
\dref{defstratum}) strata. Consequently, the union of the $l$ smallest
strata is a closed set.     

\subsection{Outline of the construction of the metric} 

The construction is based on the following observation. Let $\sigma$ be 
an element of $\Sigma_f$, and let $x_{1},\ldots,x_{m}$ be local coordinates 
defined on a \nd $U$ of $\sigma$ in $M$ \wrt which $f$ has the following 
expression~:

$$
f(x_{1},\ldots x_{m}) = x_{1}^{2} + \ldots + x_d^2 +
f(0,\ldots,0,x_{d+1},\ldots ,x_{m})\;, 
$$
with $d>0$ (\lref{tuba} and \pref{positivity}). Let us introduce
the following piece of notation~: given an open subset $O$ of $M$,~let
$$
\W_O(\Sigma_f) = \{ x \in O \; ;\; x^+ \in \Sigma_f \cap O\;\; \dim{and}\;\; 
\varphi^t(x) \in O \;\;\forall t\geq 0\} \;,
$$
where $\varphi^t$ denotes the flow of the \lw gradient vector field 
$\nabla_\F f$ of $f$ associated to the metric $g$ (\dref{lwgrdef}), and 
where $x^+ = \lim_{t\to\infty}\varphi^t(x)$ (cf.~Definition
\ref{violoncelle}). The 
set $\W_O(\Sigma_f)$ is the stable set (\dref{stable}) of $\Sigma_f
\cap O$ \wrt $(\nabla_\F f)|_{O}$. If instead of the metric $g$, we use 
the Euclidean metric \wrt the coordinates $x_{1},\ldots,x_{m}$ to 
construct $\nabla_\F f$,~then  

\begin{equation}\label{confine}
\W_U(\Sigma_f) \subset \{x_1 = 0 , \ldots , x_d =  0\} \;,
\end{equation}
and the plane $\{x_1 = 0 , \ldots , x_d =  0\}$ is transverse to
$\F$, as~needed. \\

The problem is that it is not clear that a {\it global} Riemannian
metric $g$ exists with the property that every point $\sigma$ in
$\Sigma_f$ admits local coordinates as above, \wrt which $g$ is
Euclidean (at least along the leaves). Fortunately, we do not need that
much~: for the property (\ref{confine}) to hold, it is sufficient to
have a Riemannian metric on $U$ for which the plane $\{x_1 = 0 ,\ldots
, x_d = 0\}$ is orthogonal to the planes $\{x_{d+1} = c_{d+1} ,\ldots
, x_m = c_m\}$ (cf.~\pref{Rossini}). This constitutes the key
observation, as it can be globalized, at least to a \nd of a
stratum. \sbref{mns} shows how. Briefly, starting with the (constant
rank) distribution $V^+$ on a stratum $S$, we construct a foliation
$\G\subset \F$ on a \nd $\N$ of $S$ such that $T_S\G = V^+$ (the
foliation $\G$ plays the role of the planes $\{x_{d+1} =
c_{d+1},\ldots , x_m = c_m\}$). Then, we consider the manifold $P$ of
critical points of $f$ along the leaves of $\G$ (the manifold $P$
plays the role of the plane $\{x_1 = 0,\ldots ,x_d =0\}$). It is shown
that any metric for which $TP$ is orthogonal to $T_P\G$ has the
property that $\W_\N(\Sigma_f) \subset P$ (after eventually reducing
$\N$ slightly).   

Going from one stratum to the entire $\Sigma_f$ is done by induction 
on the order of the strata (\sbref{ccsf}). Basically, to carry out this 
induction, we need to construct the foliation $\G$ in such a way that it 
``matches'' the foliations already constructed near smaller order 
strata, meaning that where two such foliations coexist, the largest 
dimensional one contains the other one. The ordering is, roughly speaking,
according to the number of strata that a given stratum contains in its
closure.  Closed strata are the smallest, strata having only closed strata
in their closure come next, and so forth. This is dictated by the
lexicographical ordering of the Boardman symbols, combined with the
rank of the distribution $V^+$. Once a ``coherent'' set of foliations
has been constructed, one may build a metric for which 
$\W_\N(\Sigma_f) \subset P$ for all strata $S$ (cf.~\sbref{cmcsf}).   

\subsection{The metric near one stratum}\label{mns}

Let $S$ be a (neither necessarily connected nor closed) embedded
submanifold of $\Sigma_f$. Suppose that $S$ intersects $\F$
transversely, and that the dimension of $V^{+}$ is constant on
$S$. Typically, $S$ is a stratum of $\Sigma_f$
(cf.~\dref{defstratum}). Consider also $\exp : \ti \OO\subset TM \to
M$, the exponential map associated to $g$, defined on a fiberwise 
convex neighborhood $\ti \OO$ of the $0$-section in $TM$. The
construction of the metric near $S$ is divided into the five
following steps.\\

\noindent
{\bf First step.} The submanifold $S$ is extended to an embedded
submanifold $S'$ of $M$, transverse to $\F$, and of dimension $q =
\,{\rm cod}\,\F$. \\

If $E = (TS + T_{S}\F)^{\perp}$, then $\exp|_{\ti \OO\cap E}$ is an 
embedding provided $\OO = \ti \OO \cap E$ is sufficiently
small. Define $S'=\exp(\OO)$. Because $\exp_{*}=\dim{Id}$ along the
$0$-section, the bundle $T_{S}S'$ coincides with $TS \oplus  E$. In
particular, the submanifold $S'$ is transverse to $\F$ along $S$. We
may assume that $\OO$ has been chosen small enough for $S'$ to be 
transverse to $\F$ everywhere.  \\ 

\noindent
{\bf Second step.} Let $\d$ be a subbundle of $T_{S}\F$ (e.g.~$\d =
V^+$ when $S$ is a stratum). We extend $\d$ to a subbundle
$\d'$ of $T_{S'}\F$. \\

Let $\nabla$ be a linear connection on $T_{S'}\F$. For every $x$ in
$S'$, there is a natural path between $x$ and a point $x_o$ in $S$.
Indeed, take $\gamma_x:[0,1]\to S':t\mapsto \exp(tX)$, where
$\exp(X)=x$. The connection $\nabla$ and the path $\gamma_x$ induce a
linear isomorphism $i_x : T_{x_o}\F\to T_x\F$, obtained by parallel
translation along $\gamma_x$ \wrt $\nabla$. For $x$ in $S'$, define
$\d'_x$ to be $i_x(\d_{x_o})$. \\

\noindent
{\bf Third step.} Given an embedded submanifold $S'$, complementary to 
$\F$, and a subbundle $\d'$ of $T_{S'}\F$, we extend $\d'$ to a smooth
foliation $\G$, defined on a \tnd $\N$ of $S'$, tangent to $\F$.\\

Consider the map 
$$
\xi : \U \subset T_{S'}\F \to M : X_s \in T_s\F \mapsto 
\tau_{X_s}(1)\, ,
$$
where $\tau_{X_s}$ denotes the geodesic for the Riemannian metric 
$g|_{F_s}$ starting at $s$ and tangent to $X_s$ (in particular 
$\tau_{X_s}\subset F_s$), and where $\U$ is some fiberwise convex \nd 
of the $0$-section on which $\xi$ is an embedding. The map $\xi$ is 
called hereafter the {\em foliated exponential map associated to the 
metric $g$}. Consider now the foliation $\ti \G$ of $\U$ by the traces of
the cosets of $\d'$, and the push-forward $\G$ of $\ti \G$ via the map 
$\xi$. Let $\N = \xi(\U)$.\\

\noindent
{\bf Fourth step.} Suppose that the foliated second differential $\ol
d{}^2\!f$ of $f$ is positive definite on $\d$. Let $P$ be the foliated
singular locus of $f$ \wrt the foliation $\G$, that is, the set $\{x\in
\N\,;\, df(T_{x}\G)=0 \}$. Provided the \nd $\N$ of $S$ is
sufficiently small, the set $P$ is an embedded submanifold of $\N$
transverse to $\G$ and closed in $\N$. \\ 

Let $p:BT\G\to\N$ be the bundle of basis of $T\G$, and define the 
map
$$
\eta : BT\G\to\R^{d} : \{e_{1},\ldots,e_{d}\} \mapsto  (df(e_{1}),
\ldots,df(e_{d}))\, ,
$$
where $d$ is the dimension of $\G$. Since $\ol d{}^{2}\!f$ is nondegenerate 
on $\d=T_S\G$, the map $\eta$ is a submersion near $p^{-1}(S)$. Therefore, 
the map $\eta$ is a submersion near $\eta^{-1}(0)$, provided $\N$ is 
sufficiently small. The submanifold $\eta^{-1}(0)$ coincides with $p^{-1} 
(P)$, where $P$ is a closed embedded submanifold of $\N$. To prove
that $P$ is transverse to $\G$, consider a vector $X$ in $T_{p}P \cap
T_{p}\G$, with $p$ in $P$. Extend $X$ to a section $\tilde X$ of
$T\G$ defined on a \nd $U$ of $p$ in $\N$. On $P\cap U$, the function
$\tilde X f$ vanishes, and because $X$ is tangent to $P$, we have $0 = 
X(\tilde X f)= \ol d{}^2\!f(X,X)$ (notice that $\ol d{}^2\!f$ is indeed
well-defined on $T_P\G$; it is the foliated second differential of $f$
\wrt the foliation $\G$). Since $\ol d{}^2\!f$ is nondegenerate 
on $T_{S}\G = \d$, it remains nondegenerate on $T_P\G$ provided $\N$
is sufficiently small. Thus, the vector field $X$ must vanish.\\

\noindent
{\bf Fifth step.} Let $h$ be a Riemannian metric for which $TP$ is
orthogonal to $T_P\G$. Such a metric can be constructed by means of a 
partition of unity for instance. The leafwise gradient vector field 
$\nabla_\F f$ of $f$ associated to $h$ (\dref{lwgrdef}) satisfies the
the following property.

\begin{prop}\label{Rossini} There exists a \nd $\M$ of $\N\cap\Sigma_f$ 
contained in $\N$ for~which 
$$
\W_{\M}(\Sigma_f)\subset P\, .
$$ 
\end{prop}

\Pf First observe that if $x$ belongs to $P$, then $(\nabla_\F f)_x$ belongs 
to $TP$ (as implied by the identity $h(\nabla_\F f, T_{x}\G) = df(x)(T_{x}\G) 
= 0$). We will then show that for some sufficiently small \nd $\M$ of $\N\cap
\Sigma_f$, no trajectory of $(\nabla_\F f)|_{\M}$ starting outside $P \cap \M$
can approach $P \cap \M \supset \Sigma_f$ in forward time. \\ 

Since $\G \subset \F$, and since $P$ is complementary to $\G$, given $x$ in
$P$, there exist coordinates $x_{1},\ldots,x_{m}$ on a \nd $U$ of $x$
such~that   
\begin{enumerate}
\item[i)] $x_{n+1}=c_{n+1},\ldots,x_{m}=c_{m}$ define $\F|_{U}$, 
\item[ii)] $x_{d+1}=c_{d+1},\ldots,x_{m}=c_{m}$ define $\G|_{U}$,  
\item[iii)] $x_{1}=0,\ldots,x_{d}=0$ define $P\cap U$,
\item[iv)] $f(x_1,\ldots,x_m)=  x_{1}^{2} + \ldots + x_{d}^{2} + f(0,
\ldots ,0, x_{d+1}, \ldots,x_{m})$. 
\end{enumerate}
The proof of this fact, very similar to that of \lref{tuba}, is
omitted.\\  

Let $x$ be an element of $P$ and let $(x_1,\ldots ,x_m)$ be nice
coordinates defined on a \nd $U$ of $x$. Let also $u= \sum^d_{i=1} u_i
\fr{\partial}{\partial x_i}$ be a nonvanishing constant vector field 
on $U$ tangent to the foliation $\G$. In the next paragraph, a dot
$\cdot$ will denote the Euclidean scalar product \wrt the coordinates
$(x_1,\ldots ,x_m)$. We have the following sequence of identities~: 
\begin{eqnarray*}
\frac{\partial}{\partial u}(\nabla_\F f\cdot u)(x) &=&  
\sum^d_{i,j=1} u_i\frac{\partial}{\partial x_{i}}
\left(\sum^n_{k=1} \frac{\partial f}{\partial x_{k}} h^{kj} u_j
\right) (x)\\  
&=& \sum^d_{i,j=1} \sum^n_{k=1} u_i u_j \left( \frac{\partial^{2}
f}{\partial x_{i}\partial x_{k}}(x) h^{kj}(x) + \frac{\partial
f}{\partial x_{k}}(x) \frac{\partial h^{kj}}{\partial x_{i}}(x) \right)  \\
&=& 2 h(u,u)(x) + \sum^d_{i,j=1} \sum^n_{k=1} u_i u_j  \frac{\partial
f}{\partial x_{k}}(x) \frac{\partial h^{kj}}{\partial x_{i}}(x)\;,
\end{eqnarray*}
where $h^{kj}$ denotes the component $k, j$ of the inverse of the matrix
of $h|_{T\F \times T\F}$. In particular, if $x$ belongs to $P\cap
\Sigma_f = \N \cap \Sigma_f$, we~have 
$$
\frac{\partial}{\partial u}(\nabla_\F f\cdot u)(x) = 2h(u,u)(x) >0\, .
$$
Hence, for all $x$ in $P\cap\Sigma_f$, there exists an open \nd $U_x$ of 
$x$, endowed with nice coordinates, such that $\frac{\partial
}{\partial u}(\nabla_\F f\cdot u)$ is strictly positive on $U_x$, for
all nonvanishing constant vector field $u$ tangent to the foliation $\G$.
As a consequence, a trajectory $\theta_x$ of $\nabla_\F f$ for which 
$\theta_x(t)$ lies in $U_x - (P \cap U_x)$ for all $t$ in $[t_o, \infty)$
cannot converge to a point in $P \cap U_x$. Indeed, Let $w$ be any
point in $U_x - (P \cap U_x)$, and let $u$ be a nonvanishing constant
vector field parallel to the line joining $w = (\sq x m)$ to
$w_o = (0, \ldots, 0, x_{d+1}, \dots, x_m)$. Then, since $\frac{\partial
}{\partial u} (\nabla_\F f\cdot u)|_{U_x} > 0$ and $(\nabla_\F
f\cdot u)(w_o) = 0$, the quantity $(\nabla_\F f\cdot u) (w)$ is
strictly positive. Supposing $w = \theta_x(t)$ for some $t\geq t_o$, the
last assertion implies that the Euclidean distance between $\theta_x(t)$
and $P \cap U_x$ grows with~$t$.\\ 

Let $\{U_i;i\geq 1\}$ be a locally finite refinement of the covering
$\{U_x ; x\in P\cap\Sigma_f\}$ of $P\cap\Sigma_f$. Let $\M = \N \cap
(\cup_{i\geq 1} U_i)$. Any trajectory of $(\nabla_\F f)|_{\M}$
starting in $P \cap \M$ remains in $P \cap \M$, and no trajectory
starting outside $P \cap \M$ will approach $P \cap 
\M$ in forward time. Indeed, suppose on the contrary that $\lim_{t\to
\infty}\theta_x(t) \in P$, for some trajectory $\theta_x$ of
$(\nabla_\F f)|_\M$. Then, for some $t_o$ and some $i$, the point
$\theta_x(t)$ lies in $U_i$ for all $t\geq t_0$. The discussion in the
previous paragraph implies that $\theta_x(t) \in P$ for all $t\geq
t_0$, and hence that $\theta_x \subset P$.   
\cqfd

\begin{rmk}{\rm The reason for introducing $\M$ instead of supposing
once more $\N$ small enough is that we need to make sure that the
portion of $\Sigma_f$ whose stable set is taken care of ($\N \cap
\Sigma_f$ in the previous proposition) is fixed once $\G$ is given.
}\end{rmk}

\begin{ccl}{\rm Suppose given the following data~:
\begin{enumerate} 
\item[-] an embedded submanifold $S$ of $\Sigma_f$ intersecting $\F$
transversely,
\item[-] a subbundle $\d$ of $T_S\F$ restricted to which $\ol d{}^2\!f$
is positive definite. 
\end{enumerate}
We can then construct a foliation $\G$ tangent to $\F$, defined on a
tubular neighborhood $\N$ of $S$, satisfying the following properties.
  \begin{enumerate}
\item[i)] $T_S\G$ coincides with $\d$.
\item[ii)] The foliated singular locus $P$ of $f$ \wrt the foliation $\G$,
           is an embedded submanifold transverse to $\G$ (equivalently
           $\ol d{}^2 \!f$ is positive definite on~$T_P\G$). In
           particular, $P$ as codimension at least one and intersects $\F$
           transversely. 
\item[iii)] Let $h$ be a Riemannian metric on $\N$ for which $TP$ is 
            perpendicular to $T_P\G$. Then there exists a \nd $\M$ of
            $\N\cap\Sigma_f$ contained in $\N$ such that the  
            leafwise gradient vector field of $f$ associated to $h$
            satisfies the property that $\W_{\M}(\Sigma_f)
            \subset P$.   
  \end{enumerate} 
}\end{ccl}

\subsection{Construction of a complete  system of
foliations}\label{ccsf} 

As explained in the beginning of \sref{metric}, the construction of
a suitable Riemannian metric goes by defining a collection of
foliations, one near each stratum, in a compatible way. The word {\it
compatible} means that where two such foliations coexist, the one 
whose dimension is largest contains the other one. This construction
is the subject of the present subsection.

\begin{dfn}\label{sf} Consider a collection of strata $S_1,\ldots,S_l$
satisfying the property that for all $k$, the set $S_k\cup \ldots \cup
S_l$ is closed in $M$. {\em A system of foliations on $\cup_{k=1}^l
S_k$} is a collection $\{ \G_k ; k=1,\ldots,l \}$ of foliations, where
$\G_k$ is defined on a \tnd $\N_k$ of an open subset $O_k = \N_k \cap
S_k$ of $S_k$. Moreover, the following properties are required to
hold~: 
\begin{enumerate}
\item[a)] For all $1\leq k\leq l,$ $S_k \subset \N_k \cup
          \ldots \cup \N_l$.
\item[b)] If $S_{k_1}$ and $S_{k_2}$ are such that $\ol S_{k_1}\cap
          S_{k_2} = \emptyset = S_{k_1}\cap \ol S_{k_2}$ then 
          $\N_{k_1}\cap\N_{k_2}=\emptyset$.  
\item[c)] ${\rm dim}\, \G_k = {\rm rk}\, (V^+|_{S_k})$.
\item[d)] If ${\rm dim}\,\G_{k_1}\geq {\rm dim}\,\G_{k_2}$, then
          $T\G_{k_1}\supset T\G_{k_2}$ on $\N_{k_1}\cap\N_{k_2}$. 
\item[e)] For all $k$, the set $P_k = \{p\in \N_k\, ;\, df(p)(T_p\G_k) = 0\}$ 
          is a submanifold of $\N_k$ transverse to $\G_k$ 
\end{enumerate} 
A {\em complete system of foliations on $\Sigma_f$} is a system of
foliations on the union of all the strata constituting $\Sigma_f$.
\end{dfn}

\begin{rmk}{\rm Notice that the set $S_1, \ldots ,S_l$ of $l$ smallest
strata, arranged in decreasing order, satisfies the property that
$S_k\cup \ldots \cup S_l$ is closed for all $1\leq k\leq l$. 
}\end{rmk}

\begin{rmk}{\rm 
Given a system of foliations on $\cup_{k=1}^l S_k$, we will
always assume (and this is not restrictive) that another system of
tubular neighborhoods $\ul{\N}_k$ of open subsets $\ul O_k =
\ul{\N}_k\cap S_k$ of $S_k$ as been given that satisfies $\ol{\ul{\N}_k}
\subset \N_k$ for all $k$, as well as property $a)$ in
\dref{sf}. Hence, the collection $\{ \G_k|_{\ul{\N}_k} \}$ is also a
system of foliations on $\cup_{k=1}^l 
S_k$. We will also use $\N'_k,$ $\N''_k, \ldots ,\N_k^{i}{}', \ldots $
to denote more tubular neighborhoods of open subsets $O'_k, O''_k,
\ldots , O_k^{i}{}', \ldots $ of $S_k$ such~that  
$$
\ol{\ul{\N}_k} \subset \N_k^{i}{}' \subset \ol{\N_k^{i}{}'} \subset
\N_k^{(i-1)}{}' \subset \N_k\;.
$$ 
The need for these additional \nds will appear in the proof of
\pref{Schubert}. It explains why, in \dref{sf}, we consider a \tnd of an
open subset of a stratum, rather than a \tnd of the entire
stratum. Indeed, if $\N$ is a \tnd of a nonclosed embedded
submanifold $S$ (like most strata), there is no \tnd $\N'$ of $S$
with $\ol \N' \subset \N$. To obtain such an inclusion we have to
replace $S$ by a submanifold $O$ of $S$ with $\ol O \subset S$.
}\end{rmk}

\begin{prop}\label{Schubert}(Recurrence step.)
Let $S_1, \ldots ,S_l$ be a set of strata satisfying the property that
for all $k$, the set $S_k\cup \ldots \cup S_l$ is closed in $M$, and
suppose that $\{\G_k;k=1,\ldots,l\}$ is a system of foliations on
$\cup_{k}S_k$. Let $S$ be another stratum for which $S \cup S_1\cup
\ldots  \cup S_l$ is closed. Then we can extend $\{ \G_k ; k=1,
\ldots, l \}$ to a system of foliations on $S\cup S_1\cup\ldots\cup
S_l$.   
\end{prop}
The word {\em extend} has to be given the following meaning. A
foliation $\G$ will be constructed on some \tnd $\N$ of an open subset
$O$ of $S$ in such a way that $\G$ together with the {\em restrictions} 
of the foliations $\G_k$ to the \nds $\ul{\N}_k$ is a  system of 
foliations on $S\cup S_1\cup \ldots\cup S_l$.\\ 

\Pf We will assume for our convenience that $\ol S\cap S_k\neq 
\emptyset$ for all $1\leq k\leq l$ (if this is not true, select the
$S_k$'s intersecting $\ol S$ nontrivially, label them $S_1, \ldots
,S_l,$ and ignore the other ones until further notice). We can also
assume that ${\rm dim}\,\G_k\geq {\rm dim}\,\G_{k+1}$ for all $1\leq
k\leq l-1$ (without affecting the property that $S_k\cup \ldots \cup
S_l$ is closed for all~$k$). \\ 

The proof follows the first three steps of \sbref{mns}. The difficulty
lies in the second and third steps. We will need to adjust the distribution 
$\d'$ and then the foliation $\G$ so as to make them match the foliations 
$\G_k$ already constructed. \\

\noindent
{\bf Adjustment of $\d'$.} Suppose $S\subset \Sigma^{(\sq i k)}_{d_o}$, 
in particular ${\rm rk}\, V^+_S = d_o$. We begin with describing a procedure 
that allows one to canonically extend a $d$-dimensional subspace of 
$T_s\F$, $s \in S$, on which $\ol d{}^2\!f$ is positive definite ($d$ may  
therefore not be greater than $d_o$), to a $d_o$-dimensional subspace with
the same property. \\ 
 
Because $\ol d{}^2\!f$ induces a metric on $T_S\F$ whose rank and signature
are constant, it yields a bundle decomposition $T_S\F = V^+\oplus
V^-\oplus V^0$ into $g$-orthogonal and $\ol d{}^2\!f$-orthogonal
subspaces (cf.~beginning of \sref{metric}). Let $p_1$ denote the
projection $T_S\F \to V^+$, and let $p_2$ denote the projection $T_S\F
\to V^-\oplus V^0.$ Given an element $X$ in $T_S\F$, write $X=X^+ +
X^- + X^0$, where $X^+$, $X^-$ and $X^0$ belong to $V^+$, $V^-$ and
$V^0\,$ respectively. For $1\leq d\leq d_o$, consider   
\begin{eqnarray*}
\W^+ & = & \cup_{s\in S} \left(\W^+_s = \{X\in T_s\F\, ;\,
\ol d{}^2\!f(X,X) > 0\}\right) \\
\p^{+, d} & = &  \cup_{s\in S} \left(\p^{+, d}_s = \{P\in
G^d(T_s\F)\, ;\, P-\{0\}\subset\W^+_s\}\right)\;,
\end{eqnarray*}
where $G^d(T_s\F)$ denotes the Grassmann manifold of $d$-planes in
$T_s\F$. If $E$ is a vector space endowed with an inner product, let
$S(E)$ denotes the unit sphere in $E$, and let $B(E)$ denotes the open
unit ball in $E$. Observe that if $P$ is a plane in $\p^{+, d}$ then
$p_1|_P$ is injective. Hence, an element $P$ of $\p^{+, d}_s$ is the
graph of a linear map defined on~$p_1(P)$~:
$$\begin{array}{ccccl}
\varphi_P & : & p_1(P) & \to & V^-_s\times V^0_s \\
          &   &   x   & \mapsto & p_2\circ (p_1|_P)^{-1}(x)\, ,      
\end{array}$$
whose restriction to the sphere $S(p_1(P))$ takes its values in
$B(V^-_s) \times V^0_s$. Conversely, an element of $G^d(T_s\F)$ that
coincides with the graph of such a linear map belongs to $\p^{+,
d}_s$. Notice that if $P$ is $d_o$-dimensional, then $p_1(P) = V^+_s$.\\ 

Now given an element $P$ in some $\p^{+,d}_s$ with $d<d_o,$ we extend it
to an element $\e(P)$ in $\p^{+,d_o}_s$ as described hereafter.
Let $q_P : V^+_s \to p_1(P)$ denote the orthogonal projection
onto $p_1(P).$ Then define $\e(P)$ via its associated linear map
$\varphi_{\e(P)}$~by~:
$$
\begin{array}{ccccl}
\varphi_{\e(P)} & : & V^+_s & \to & V^-_s \times V^0_s \\
                &   & x     & \mapsto & (\varphi_P\circ q_P)(x)\;. 
\end{array}
$$

\begin{rmk}\label{clc}{\rm We can form convex linear combinations of
elements in $\p^{+,d_0}$. Indeed, let $P_1, \ldots, P_r$ be elements of
$\p^{+,d_0}_s$ for some $s$ in $S$, and let $a_1,\ldots ,a_r$ be positive
real numbers with $a_1 + \ldots + a_r = 1$. Then $a_1 P_1 + \ldots + 
a_r P_r$ denotes the element of $\p^{+,d_0}_s$ whose associated function 
$\varphi$ is $a_1 \varphi_{P_1} + \ldots + a_r\varphi_{P_r}$. The map
$\varphi$ takes its values in $B(V^-_s) \times V^0_s$, as this set is
convex.  
}\end{rmk}

As in the first step of \sref{mns}, the stratum $S$ is extended to a
$q$-di\-men\-sio\-nal embedded submanifold $S'$ transverse to $\F$. The
manifold $S'$ is the image of a fiberwise convex \nd $\OO$ of the
$0$-section in the bundle $p : E = (TS + T_S\F)^{\perp} \to S$ via the
exponential map associated to $g$. Recall from the second step
of \sref{mns}, that a linear connection on $T\F$ determines, for each
point $x$ in $S'$, a linear isomorphism $i_x : T_{x_o}\F\to T_x\F$,
where $x_0 = \exp(p(\exp^{-1}(x)))$. Set $\OO'_k = \exp (\OO\cap
p^{-1}(\N'_k \cap S))$ and $\OO''_k = \exp (\OO \cap p^{-1}(\N''_k
\cap S))$.

\begin{lem}\label{flute} Provided $\OO$ is sufficiently small, the
following properties hold for all~$k$~:
\begin{enumerate}
\item[-] $\OO'_k \subset\N_k$, 
\item[-] $\OO''_k\supset \N'''_k \cap S'$, 
\item[-] for every $x$ in $\OO'_k$, the space $i_x^{-1}(T_x\G_k)$ is
contained in $\W^+_{x_o} \cup \{0\}$.  
\end{enumerate}
\end{lem}

\Pf Given a subset $U$ of $S$ and a positive number $\varepsilon$, let
$U^\varepsilon$ denote the subset $\exp (\cup_{s\in U} B_{\varepsilon}
(0_s))$ of $S'$, where $B_\varepsilon(0_s)$ is the ball of radius 
$\varepsilon$ centered at $0$ in $p^{-1}(s)$. For every $s$ in $S$,
there exists a \nd $U_s$ of $s$ in $S$ and a positive number
$\varepsilon_s$ such that 
\begin{enumerate}
\item[-] $U^{\varepsilon_s}_s\subset \N_k$ for all $k$
         for which $s \in \ol{\N}'_k$, 
\item[-] $U_s \cap \N'_k = \emptyset$ for all $k$ for
         which $s\notin\ol{\N}{}'_k$,
\item[-] $U_s\subset \N''_k$ for all $k$ for which $s\in\ol{\N}'''_k$,
\item[-] $U^{\varepsilon_s}_s \cap \N'''_k = \emptyset$ for all $k$ for
         which $s \notin \ol{\N}'''_k$,    
\item[-] $i_x^{-1}(T_x\G_k)\subset \W^+_{x_0} \cup \{0\}$ for all $x\in
U^{\varepsilon_s}_s$ with $s$ in $\ol{\N}{}'_k$.  
\end{enumerate}
Let $\{U_i;i\geq 1\}$ be a locally finite refinement of the covering
of $S$ by the $U^{\varepsilon_s}_s$'s. For any fiberwise convex \nd
$\OO$ of the $0$-section in $E$ such that $\exp(\OO) \subset \cup_iU_i$,
the required three properties are satisfied.  
\cqfd  

\vspace{.2cm}
For all $k=1,\ldots,l$, the dimension of $\G_k$ is at most equal to
$d_o$ (we supposed that $\ol S\cap S_k  \neq \emptyset$ for all $k$,
and this implies that $\dim{rk} V^+_S \geq \dim{rk}V^+_{S_k}$ for all
$k$). We can therefore associate to each foliation $\G_k$ a
$d_o$-dimensional distribution $\d'_k$, defined on $\OO'_k$,
tangent to $\F$, and containing $T\G_k$. Simply let 
$$
(\d'_k)_x = i_x \left( \e\left(i_x^{-1} (T_x \G_k) \right)\right)\,.  
$$ 

We would like to paste the $\d'_k$'s together so as to obtain a
$d_o$-dimensional distribution $\d' \subset T_S'\F$ on $S'$ containing
the distribution $T_{S' \cap \N'''_k}\G_k$, and for which $\ol
d{}^2\!f$ is positive definite on $\d = \d'|_S$. For each $k$ between
$1$ and $l$, let $\rho_k : M\to [0,1]$ be a smooth function equal to
$1$ near $\N''_k$ and vanishing near the complement of $\N'_k$. We
define $\d'$ as~follows~:   

$$
\begin{array}{ccl}
\d'_x & = & \rho_1 (x_o)(\d'_1)_x + (1-\rho_1(x_o))\rho_2(x_o)(\d'_2)_x
                   + \ldots + \\
             &   & (1-\rho_1(x_o)) \ldots (1-\rho_{k-1}(x_o)) 
                   \rho_k(x_o) (\d'_k)_x + \ldots + \\
             &   & (1-\rho_1(x_o)) \ldots (1-\rho_l(x_o))
                   i_x (V^+_{x_o}) \;,
\end{array}
$$
where $x$ is in $S'$, and where $x_o = \exp (p(\exp ^{-1}x))$. It is understood 
that when some distribution $\d'_k$ is not defined at the point $x$, the 
quantity $\rho_k(x_o) (\d'_k)_x$ is defined to be $\rho_k(x_o) i_x 
(V^+_{x_o})$. The linear combinations appearing in the right hand side 
has to be understood as follows. If $P_1, \ldots, P_r$ are subspaces of 
$T_x\F$ for some $x$ in $S'$ such that $i_x^{-1} (P_i) \subset \W^+ \cup 
\{0\}$, and if $a_1,\ldots ,a_r$ are real numbers, then    
$$
a_1 P_1 + \ldots + a_r P_r = i_x \left( a_1 i_x^{-1} (P_1) + \ldots
+ a_ri_x^{-1} (P_r)\right)\, ,
$$ 
where the linear combination in the right hand side has to be interpreted 
according to \rref{clc}. Observe that on $\OO''_k$, hence on $\N'''_k\cap 
S'$, the distribution $\d'$ contains $T\G_k$. Indeed, the function $(1- 
\rho_{k}(x_o))$ vanishes on $\N''_k$, hence, only $\d'_1, \ldots ,\d'_k$ 
are involved in the definition of $\d'_x$ for $x$ in $\OO''_k$. By construction, 
$\d'_j$ contains $T\G_k$ for $j$ between $1$ and $k$; hence, any convex 
linear combination of $\d'_1, \ldots ,\d'_k$ contains $T\G_k$ as~well. \\   

To keep notations light, we denote $\N'''_k$ by $\N_k$, $\N''''_k$ by
$\N'_k$, and so forth, while fixing~$\ul{\N}_k$.\\ 

\noindent
{\bf Adjustment of $\G$.} We will now extend the distribution $\d'$ to 
a foliation $\G \subset\F$ defined on a \tnd $\N$ of $S'$, that
contains $\G_k$ on $\N \cap \N'''_k$. Recall from the previous section 
the foliated exponential map~:    
$$
\xi : \U \subset T_{S'}\F \to M\;,
$$
where the set $\U$ is a fiberwise convex \nd of the $0$-section in
$T_{S'}\F$, small enough for $\xi$ to be an embedding. Let
$\tilde{\G}_k$ denote the foliation $(\xi)^{-1}_*(\G_k)$ defined on 
$\xi^{-1}(\N_k)$. It is of course tangent to the foliation of $\U$ by
the fibers of the natural projection $\pi : T_{S'}\F\to S'$. For all 
$k=1,\ldots, l$, let $\e_k = (T\G_k)^{\perp} \cap T\F$. Given any $k$, 
there exists a (not necessarily fiberwise convex) \nd $\V_k$ of the 
$0$-section of $T_{S'\cap\N_k}\F$ contained in $\xi^{-1}(\N_k)$, 
such that each leaf of $\tilde{\G}_k|_{\V_k}$ intersects $\e_k$ along exactly 
one point. In particular, there is a map  $g_k: \V_k \to \e_k$ defined by 
$\{g_k(x)\} = (\tilde G_k)_x\cap \e_k$, where $(\tilde G_k)_x$ denotes 
the leaf of $\tilde{ \G}_k|_{ \V_k}$ containing $x$. The following
lemma is very similar to \lref{flute} although it is formulated in
$T_{S'}\F$ instead of $M$. Its proof is omitted. Let $\U'_k$
(respectively $\U''_k$) denote the set $\U \cap \pi^{-1}(\N'_k \cap
S')$ (respectively $\U \cap \pi^{-1}(\N''_k \cap S')$). 
 
\begin{lem} Provided $\U$ is sufficiently small, the following
properties hold for all $k$. 
  \begin{enumerate}
\item[-] $\U'_k \subset \V_k$, 
\item[-] $ \U''_k \supset \xi^{-1}(\N'''_k)$.
  \end{enumerate}
\end{lem}

Let $g_{l+1}$ denote the orthogonal projection $T_{S'}\F \to \d'^{\perp}$, 
and let $\nu : T_{S'}\F \times T_{S'}\F\to T_{S'}\F$ be the map that sends 
the pair $(x,y)$ to the orthogonal projection of $x$ onto the coset of $\d'$
passing through the point $y$ (that is, $\nu(x,y)=x + g_{l+1} (y-x)$). 
Consider a smooth function $\rho_k : M\to [0,1]$ whose value is $1$ near 
$\N''_{k}$, and is $0$ near the complement of $\N'_{k}$. Define $\mu : 
\U \to T_{S'}\F$~by     

$$\begin{array}{ccl}
\mu(x) & = & \nu \left(x\;,\;\rho_1 (x_o) g_1(x) + (1-\rho_1 (x_o)) \rho_2
             (x_o) g_2 (x) + \ldots +\right. \\
       &   & (1 - \rho_1(x_o)) \ldots (1 - \rho_{k-1} (x_o))\rho_k(x_o)
             g_k(x) + \ldots + \\ 
       &   & \left.(1 - \rho_1(x_o)) \ldots  (1 - \rho_l
             (x_o))g_{l+1}(x)\right) \;, 
\end{array}$$
where $x_o=\pi(x),$ and where, if $g_k(x)$ is not defined, we set 
$\rho_k (x_o) g_k(x) = 0_{x_o}$. Observe that the map $\mu$ coincides with 
the identity on the $0$-section. Denote by $\G'$ the foliation of 
$T_{S'}\F$ by the cosets of $\d'$. The idea is that, provided $\mu$ is
a diffeomorphism, the foliation $\mu (\ti{\G}_k|_{\U''_k})$ is tangent
to the foliation $\G'$. 

\begin{lem} There exists a fiberwise convex \nd $\ti\U \subset \U$ of
the $0$-section in $T_{S'}\F$ for which $\mu : \ti \U \to T_{S'}\F$ is
a diffeomorphism onto its image, another \nd of the $0$-section in
$T_{S'}\F$.  
\end{lem}

\Pf
First observe that the map $\mu_{k} (x) = \nu(x,g_k(x))$, defined for 
$x$ in $\V_k$, coincides with the identity on $\V_k \cap \e_k$ (or on
$T_{S'}\F$ for $k=l+1$), and that $(\mu_k)_*(X) = X$ for $X$ in
$T_s \G_k$, with $s$ in $S'\cap \V_k$ (as a consequence of the
fact that $(\d')_s$ contains $T_s\G_k$ for $s$ in $S' \cap
\N_k$). Hence, the map $(\mu_k)_{*_s}$ coincides with the identity map
for $s$ in $S' \cap \V_k$. Observe now that, if we fix $s$ in $S'$ and
let $x$ vary in $T_s\F\cap\U$, the map $\mu$ can be written as~:   
$$
\mu(x) = \nu(x,t_1 g_1 (x) + \ldots + t_l g_l (x) + t_{l+1} g_{l+1}(x))\; ,
$$
where $t_k = (1-\rho_1(s))\ldots (1-\rho_{k-1}(s))\rho_k(s)$ for $1
\leq k \leq l$, and where $t_{l+1}=(1-\rho_1(s))\ldots (1 - \rho_l
(s))$. Hence, for $X$ in $T_{0_s}T_s\F \simeq T_s\F$, we have   
\begin{eqnarray*}
\mu_{*_s} (X) & = & \nu_{*_{(s,s)}}\left(X,\left(t_1g_1+\ldots + 
              t_{l+1}g_{l+1}\right)_{*_s}(X)\right)\\
              & = & \nu_{*_{(s,s)}}\left((t_1 + \ldots + t_{l+1})\, 
              X, t_1 (g_1)_{*_s}(X) + \ldots + t_{l+1} (g_{l+1})_{*_s}(X) 
              \right)\\
              & = & t_1 \;\nu_{*_{(s,s)}}\left( X,(g_1)_{*_s}(X)\right) + 
              \ldots + t_{l+1} \;\nu_{*_{(s,s)}}\left(
              X,(g_{l+1})_{*_s}(X) \right)\\ 
              & = & t_1 \;(\mu_1)_{*_s}(X) + \ldots + t_{l+1}\;
              (\mu_{l+1})_{*_s} (X)\\
              & = & t_1 X + \ldots + t_{l+1} X\\
              & = & X\;.
\end{eqnarray*}
Thus $\mu_{*_{s}} = \dim{Id}$ for all $s$ in $S'$, and the lemma follows.
\cqfd

\begin{lem}\label{Haydn} For all $k$, the foliation $\mu_*^{-1}(\G')$
contains the foliation $\tilde{\G}_k$ on $\U''_k\cap \ti \U$, hence on
$\xi^{-1}(\N'''_k) \cap \ti \U$.  
\end{lem}

\Pf If $x$ belongs to $\U''_k\cap\ti \U$, then $(1-\rho_k(x_o))=0$. 
Hence the map $\mu$ can be expressed as follows~:
\begin{eqnarray*}
\mu (x) & = & \nu \left( x,\;\rho_1 (x_o) g_1(x) + (1-\rho_1 (x_o)) \rho_2
(x_o) g_2 (x) + \ldots\; +\right.\\
        &   & \left.(1-\rho_1 (x_o)) \, \ldots \,(1 - \rho_{k-1}
(x_o)) \rho_k (x_o) g_k(x)\right) \, .
\end{eqnarray*}
Besides, for all $j=1,\dots, k$, the map $x\mapsto\rho_j(x_o) g_j(x)$ is
constant on the leaf $(\ti G_k)_x$ of $\tilde{\G}_k|_{\U''_k\cap\ti
\U}$ passing through $x$. This assertion follows from the fact that
whenever $k_1>k_2$, the distribution $T\G_{k_1}$ contains the
distribution $T\G_{k_2}$ on $\N_{k_1} \cap \N_{k_2}$. Thus the map
$\mu$ sends the entire leaf $(\ti G_k)_x$ into a coset of $\d'$.    
\cqfd 

\vspace{.2cm}
The foliation $\G$ is defined to be $\xi_*(\mu^{-1}_* (\G')|_{\ti \U})$.
Its domain is the open set $\N=\xi(\ti\U)$. \lref{Haydn} implies that
$\G$ contains $\G_k$ on $\N'''_k \cap \N$. Define $P$ to be the
foliated singular locus of $f$ \wrt $\G$, that is, $P=\{p \in \N ;
df(p) (T_p\G) \}$. We may assume, after shrinking $\N$ if needed, that
$\ol d{}^2\!f$ is positive definite on $T_P\G$ (cf.~Fourth step of
\sbref{mns}). The strata we might have ignored in the very beginning
of this proof should now be re-incorporated.  

\begin{lem}The foliations $\G,\G_1,\ldots,\G_l$, defined on $\N,
\ul{\N}_1, \ldots , \ul{\N}_l$ respectively, form a system of foliations on
$S \cup S_1 \cup \ldots \cup S_l$.   
\end{lem}

\Pf Among the five defining properties of a system of foliations,
properties a), c), d) and e) have been taken care of during the
construction. Only b) requires some attention. We need to make sure
that $\N$ does not intersect the \tnds of the strata that we discarded
at the beginning of the proof. A way to insure this is to fix a family
of \nds $\ti \N_k$, one for each strata $S_k$, satisfying the
following property~: 

\begin{quote}
if $\ol S_{k_1}\cap S_{k_2}= \emptyset = S_{k_1}\cap \ol S_{k_2}$,
then $\ti \N_{k_1}\cap\ti \N_{k_2}=\emptyset$ as well. 
\end{quote}
Then, whenever we consider a \nd of some strata $S_k$, we
request that it be contained in~$\ti \N_k$. 
\cqfd
\noindent
{\em End of the proof of \pref{Schubert}.}

\begin{cor}\label{Webern} A complete system of foliations on
$\Sigma_f$ exists.
\end{cor}
\Pf Let $S_1,\ldots,S_l$ be the set of all strata of $\Sigma_f$,
presented in such a way that $S_1>S_2>\ldots.$ Then, as described
in \sbref{mns} (or in \pref{Schubert} with $l=0$), we can construct a
foliation $\G_l$ near $S_l$ and use \pref{Schubert} repeatedly to
extend $\G_l$ to a system of foliations near $\cup_kS_k$.   
\cqfd

\subsection{Construction of a metric from a complete system
of foliations}\label{cmcsf} 

Let $\{\G_k;k=1,\ldots,l\}$ be a complete system of foliations on
$\Sigma_f$. Each $\G_k$ is defined on a tubular \nd $\N_k$ of an open
subset $O_k$ of the stratum $S_k$. The foliated singular locus of $f$
\wrt $\G_k$ is denoted by $P_k$. We suppose that $\dim{dim}\G_k \geq
\dim{dim} \G_{k+1}$ for all $k$.      

\begin{dfn}\label{masf} A Riemannian metric is said to be {\em
adapted to the system of foliations $\{\G_k;k=1,\ldots,l\}$} if
$T_xP_k$ is perpendicular to $T_x \G_k$ for all $x$ in $P_k \cap
\N'_k$.  
\end{dfn}

\begin{prop}\label{Purcell} A Riemannian metric adapted to the system
of foliations $\{\G_k;k=1,\ldots,l\}$ exists.
\end{prop}

\Pf For $x$ in $M$, let $U_x$ be a \nd of $x$ in $M$, and let $h_x$ be
a Riemannian metric defined on $U_x$ such that 
\begin{enumerate}
\item $U_x \subset \N_k$ for all $k$ for which $x \in P_k \cap \ol{\N}{}'_k$,
\item $U_x \cap P_k\cap\ol{\N}{}'_k = \emptyset$ for all $k$ for which
      $x \notin P_k \cap \ol{\N}{}'_k$,
\item if $x \in P_k\cap\ol{\N}{}'_k$, then $T_yP_k \perp_{h_x} 
      T_y\G_k$ for all $y$ in $U_x \cap P_k$. 
\end{enumerate}
Existence of $U_x$ and $h_x$ is easily seen, except perhaps for the
last condition. Suppose that $x$ belongs to $\ol{\N}{}'_k \cap P_k$
\Iff $k \in \{k_1, k_2, \ldots,k_r\}$, with $k_1 < k_2 < \ldots < k_r$. 
Then for Property 3.~to hold, it is sufficient that the following bundles 
be pairwise $h_x$-orthogonal.    
$$
TP_{k_1}, TP_{k_2}\cap T\G_{k_1}, TP_{k_3}\cap T\G_{k_2},
\ldots, TP_{k_r}\cap T\G_{k_{r-1}}, T\G_{k_r}\;.
$$ 
These bundles span $TM$ on $U_x \cap (\cap_j(\ol{\N}{}'_{k_j} \cap
P_{k_j}))$ and are linearly independent, so that the Gram-Schmidt
orthogonalization process can be carried out. \\

let $\{U_i\, ;\, i\geq 1\}$ be a locally finite refinement of the covering 
of $M$ by the $U_x$'s, and let $\{\theta_i\}$ be a partition of unity
subordinate to the covering $\{U_i\}$. Denote by $h_i$ the metric
$h_{x_i},$ where $x_i$ is chosen in such a way that $U_i \subset 
U_{x_i}.$ Define
$$
h = \sum_{i\geq 1} \theta_i \, h_i\, .
$$
Then $T_yP_k\perp_h T_y\G_k$ whenever $y$ belongs to $P_k \cap \N'_k$. 
Indeed, let $y\in P_k \cap \N'_k$, and suppose that $y$ belongs to
$U_i \subset U_{x_i}$ for some $i$. Then $x_i\in P_k \cap 
\ol{\N}{}'_k$. Hence $T_yP_k$ is perpendicular to $T_y\G_k$ with respect 
to the metric $h_i$. Since this holds true for every $i$ for which $y\in 
U_i$, we have $T_yP_k \perp_h T_y\G_k$.   

\subsection{Conclusion}

The preceding subsections show how to construct a complete 
system of foliations $\{\G_k;k=1,\ldots,l\}$ on $\Sigma_f$ (\dref{sf}
and \cref{Webern}), as well as a Riemannian metric $h$ adapted to that
system (\dref{masf} and \pref{Purcell}). By \pref{Rossini}, there
exists, for every $k$, a \nd $\M_k$ of $\N'_k\cap\Sigma_f$ for which    
$$
\W_{\M_k}(\Sigma_f)\subset P_k\,,
$$
where $P_k$ denotes the foliated singular locus of $f$ \wrt the
foliation $\G_k$. The submanifold $P_k$ has codimension at least one
and intersect $\F$ transversely. Notice that the $\M_k$'s
cover~$\Sigma_f$.  

\section{The proof}\label{proof}

Let us recall the statement of \tref{ht} whose proof will be completed
in the present section.   

\begin{thm}\label{ht2}
On an open foliated manifold, any open relation
invariant under foliated isotopies satisfies the parametric
h-principle. 
\end{thm}

\Pf Let $\Omega$ be an open relation invariant under foliated
isotopies on the open foliated manifold $(M,\F)$. Consider a proper
$\F$-generic function $f : M\to [0,\infty)$, without leafwise local
maxima. As observed in \rref{Verdi}, we may assume, \Wloge, that $f$
is not only $\F$-generic but also strongly $\F$-generic. Let $\Sigma_f
= S_1 \cup \ldots \cup S_l$ be the decomposition of the foliated
singular locus of $f$ (\dref{lwcr}) into strata
(\dref{defstratum}). Let $\{\G_k; k=1, \ldots, l\}$ be a complete
system of foliations (\dref{sf} and \cref{Webern}), endowed with an
adapted metric $h$ (\dref{masf} and \pref{Purcell}). As before, we
will denote by $\N_k$ the domain of definition of the foliation
$\G_k$, by $P_k$ the foliated singular locus of $f$ \wrt the foliation
$\G_k$, and by $\M_k$ a \nd of $\N_k\cap \Sigma_f$ for which
$\W_{\M_k} (\Sigma_f) \subset P_k$. The \lw gradient vector field of
$f$ \wrt $h$ (\dref{lwgrdef}) is denoted as before by $\nabla_\F f$,
and its local flow by $\varphi^t$. We will need a partition
$\p=\{\,a_0=0=a_1< a_2 <\ldots<a_i<\ldots\,\}$ of $[0,\infty)$ by non
critical values (except for $a_0$ and $a_1$) of $f$. The latter
provides us with an exhaustion $M=\cup_i K_i$ of $M$ by compact
subsets $K_i=f^{-1}([a_i,a_{i+1}])$. \\   

The proof of \tref{ht2} consists in showing that, provided the
partition $\p$ is fine enough, the hypotheses of \pref{exhaustion} are
satisfied. The first hypothesis,
that the h-principle is valid near $K_0$ is easy to handle. Indeed,
$K_0 = f^{-1}(0)$ is a finite union of points (we assume here that
the minimum value of $f$ is $0$), and the local h-principle
(\dref{localhp}) is valid for any open relation (cf.~\pref{Schumann}),
in particular for $\Omega$. Concerning the second hypothesis, that the
h-principle for extensions is valid on each pair $(K_{i+1},K_i)$,
it is proven in two steps. The combination of \lref{alto} and
\lref{violon} proves the h-principle for extensions on each pair
$(K_{i+1},K_i\cup \W^{a_{i+2}}_{a_{i+1}} (\Sigma_{f}))$. It remains
to prove that the h-principle for extensions is valid on each pair
$(K_i \cup \W^{a_{i+2}}_{a_{i+1}} (\Sigma_{f}), K_i)$. This is where the
construction of a Riemannian metric carried out in \sref{metric} is
needed, as will become clear below. Fix a slice $K_i$.

\begin{ob}\label{piano}{\rm Since $f$ is a strongly $\F$-generic
function on $(M,\F)$, for each leaf $F$ of the foliation $\F$, the
critical points of $f|_F$ are isolated in $F$. Moreover, the
submanifold $\Sigma_f$ is closed and embedded. Hence, for each
$\sigma$ in $\Sigma_f$, there exists a \nd $U_\sigma$ of $\sigma$
in $M$ satisfying the following properties~:   
\begin{enumerate}
\item $U_\sigma$ is the domain of a chart $(U_\sigma,\varphi_\sigma)$
      adapted to $\F$ and centered at $\sigma$ such that
      $\varphi_\sigma( U_\sigma) = B_\sigma^1\times B_\sigma^2$, where
      $B_\sigma^1$ is a closed ball about $\varphi(\sigma)$ in the
      image of the leaf $F_\sigma$, and where $B_\sigma^2$ is a closed
      ball about $\varphi(\sigma)$ in the transverse direction.
\item $U_\sigma$ is contained in $\M_k$ for some $k$.
\item The compact set $b U_\sigma \stackrel{\dim{def.}}{=}
      \varphi_\sigma^{-1}(\partial 
      B_\sigma^1 \times B_\sigma^2)$ does not intersect $\Sigma_f$. 
\end{enumerate}
}\end{ob}
Let  $\{U_r\, ;\, r=1,2, \ldots\}$ be a locally finite refinement of the
cover of $\Sigma_f$ by the $U_\sigma$'s. Observe that $\cup_rbU_r$ is
a closed set, and that for only finitely many $r$'s, the set $U_r$
intersects $K_i$ 
nontrivially. For each $r$, choose a $\sigma_r$ in $\Sigma_f$ such that
$U_r\subset U_{\sigma_r}$. Choose also an index $k_r \in \{1, \ldots, l\}$ 
for which $U_r\subset\M_{k_r}$. We will use the notations introduced 
in \sref{lwgr}.

\begin{lem}\label{trompette} 
$$
\W(\Sigma_f)\subset \bigcup_r \; \ES (U_r \cap P_{k_r})\,.
$$ 
\end{lem}  

\Pf Let $x$ be a point in $M-\Sigma_f$ such that $x^+$ is in
$\Sigma_f$. The critical point $x^+$ belongs to $U_r$ for some
$r$. Either the trajectory $\varphi^{[0,\infty)}(x)$ intersects
$bU_r$, or it is entirely contained in $U_r$. In the second case, we
deduce from \pref{Rossini}, and from the fact that $U_r \subset
\M_{k_r}$, that $x$ must be contained in $U_r\cap P_{k_r}$. In the
first case, let $x_0=\varphi^t(x)$ be the point in
$\varphi^{[0,\infty)}(x)\cap bU_r$ for which $t$ is maximum. Then
$\varphi^{[t,\infty)}(x)$ is entirely contained in $U_r$, and
\pref{Rossini} implies that $x_0$ must belong to $P_{k_r}$. In both
cases, $x$ belongs to $\ES (U_r\cap P_{k_r})$.   
\cqfd

\begin{lem}\label{slicing}Let
\begin{eqnarray*}
\varepsilon_i & = & \dim{inf} \,\{\,f(x)-f(x^{-}) \,;\, x \in \cup_r
(bU_r\cap P_{k_r})\cap K_i \mbox{ and } x^- \in K_i\,\}\,.    
\end{eqnarray*}
Then $\varepsilon_i$ is positive.
\end{lem}

\Pf
Suppose on the contrary that $\varepsilon_i = 0$. Then there exists a 
sequence $(x_j)$ of points in $\cup_r b U_r\cap K_i$ with $x_j^-$ in
$K_i$ such that $\lim_{j \to \infty} \left(f(x_j)-f(x_j^-)\right) =
0$. We may assume, after extracting a subsequence if necessary, that  
the sequence $(x_j)$ converges in $K_i$ to a point $x$. Since $x$
belongs necessarily to $\cup_r b U_r$, it is not a leafwise critical
point of $f$, and we may therefore consider a chart $(U,\psi)$ about
$x$, adapted to the foliation $\h$ of  $M-\Sigma_f$ by the orbits of
$\nabla_\F f$. For all sufficiently large $j$, the point $x_j$ lies in
$U$. Letting $x'=\varphi^t(x)$ with $t<0$ be an element in
$U$, there exists a sequence $x'_j=\varphi^{t_j}(x_j)$ with $t_j <0$
converging to $x'$. Because $f$ is strictly increasing along
nonconstant trajectories of $\nabla_\F f$,  
\begin{eqnarray*}
0 & < & f(x) - f(x') \\
  & =  & \lim_{j\to\infty}\left(f(x_j) - f(x'_j)\right)\\ 
  & \leq & \lim_{j\to\infty} \left(f(x_j) - f(x^-_j)\right) \, ,
\end{eqnarray*}
contradicting the hypothesis that $\lim_{j \to \infty} \left( f(x_j) -
f(x_j^-) \right) = 0.$
\cqfd

Now let $\p'=\{b_0 = 0 = b_1< b_2 < \ldots<b_j\ldots\}$ be a
refinement of the partition $\p$ such that, if $a_i = b_{j_i}$, then
$b_{j+1}-b_j<\varepsilon_i$, whenever $j_i\leq j<j_{i+1}$. Fix $j$, let
$L_j$ denote the slice $f^{-1}([b_j,b_{j+1}])$, and let $i$ be the
index for which $[b_j,b_{j+1}]\subset [a_i, a_{i+1}]$. 

\begin{lem}\label{clarinette} For all $r$, the set $N_r = \ES(U_r\cap
P_{k_r} \cap L_j)\cap L_j$ is a finite union of compact subsets of
embedded submanifolds of codimension at least one, transverse to $\F$.    
\end{lem} 

\begin{rmk}{\rm It is not true in general that $N_r$ is an embedded 
submanifold. The problem is that, even though $\nabla_\F f$ is a leafwise
{\em gradient} vector field, there might be trajectories leaving $U_r\cap
P_{k_r}$ and coming back later to $U_r\cap P_{k_r}$, creating in the
process self-intersections in~$N_r$.}
\end{rmk}

{\em Proof of \lref{clarinette}.} As in the proof of \lref{slicing},
let $\h$ denote the foliation of $M-\Sigma_f$ by the orbits of
$\nabla_\F f$. For any point $x$ in $bU_r\cap P_{k_r} \cap L_j$, there
exists a $t_x \leq 0$ for which $\varphi^{t_x} (x) = x_0$ belongs to
$f^{-1}(b_j) - \Sigma_f$. Indeed, $b_{j+1}-b_j < \varepsilon_{i}$, but 
$f(x)-f(x^-)\geq \varepsilon_{i}$, by definition of $\eps_i$. Let $V_x$   
be an open \nd of $\varphi^{[t_x,0]} (x)$ in $M$ such that  
\begin{enumerate}
\item[-] $\h|_{V_x}$ is isomorphic to a product foliation, 
\item[-] each leaf of $\h|_{V_x}$ intersects $f^{-1}(b_j)$ nontrivially. 
\end{enumerate} 
Let also $U_x$ be a \nd of $x$ in $V_x\cap \M_{k_r}$ for which
$U_x\cap P_{k_r}\simeq U^1_x \times U^2_x$, where 
\begin{enumerate}
\item[-] $U^1_x$ is a \nd of $x$ in the leaf $H_x$ of $\h$ through $x$, 
\item[-] $U^2_x$ is an submanifold containing $x$ whose projection
         onto the leaf space of $\h|_{V_x}$ is an embedding.  
\end{enumerate} 
Existence of $U_x$ is guaranteed by the fact that $P_{k_r}-\Sigma_f$ is 
tangent to $\h$. Then $\ES (U_x\cap P_{k_r})\cap V_x$ is an embedded 
submanifold (it is isomorphic to $U^2_x\times \R$). Now, since $bU_r\cap 
P_{k_r} \cap L_j$ is compact, it is covered by a finite number of
$U_x$'s, say by $U_{x_1},\ldots,U_{x_n}$. Moreover, we may assume that
for each $1\leq \ell \leq n$, there exists a relatively compact
refinement $V_\ell$ of $U_{x_\ell}$, such that $\cup_\ell V_\ell$ also
cover $bU_r\cap P_{k_r} \cap L_j$. Then the compact sets $\ES
( \ol{V_\ell} \cap P_{k_r}) \cap L_j$, $\ell = 1,\ldots,n$ and
$U_r\cap P_{k_r}$ cover $N_r$, and are contained in the submanifolds
$\ES (U_{x_\ell}\cap P_{k_r}) \cap V_{x_\ell}$, $\ell = 1,\ldots,n$ and
$P_{k_r}$ respectively. Because $P_{k_r}$ is transverse to $\F$, and
because $\h$ is tangent to $\F$, the submanifold $\ES (U_{x_\ell}\cap
P_{k_r}) \cap V_{x_\ell}$ is transverse to $\F$ as well (since $U_x^2$ is
transverse to $\F$). Moreover, by construction, the submanifolds
$\ES(U_{x_\ell}\cap P_{k_r}) \cap V_{x_\ell}$ and $P_{k_r}$ have
codimension at least one.   
\cqfd 

\vspace{.2cm}
It follows from \lref{clarinette} and \lref{trompette} that the stable
set of $L_j\cap \Sigma_f$ \wrt $(\nabla_\F f)|_{L_j}$ is contained in a
finite union of compact subsets of embedded submanifolds intersecting
$\F$ transversely, and having their codimension at least equal
to one. We can thus use \tref{Chopin} and \rref{orgue} to conclude.\\

\noindent{\em End of the proof of \tref{ht2}.}

\section{Examples of open foliated manifolds}\label{ex}

The very first class of examples of open foliated manifolds consists
of the products $(M,\F) \times \R$. Let $f$ be any positive, proper,
$\F$-generic function on $(M,\F)$ (such a function always exists). The
function $g : M \times \R \to \R : (x, t) \mapsto f(x) + t^2$ has no
leafwise local maxima, and satisfies the hypotheses of \dref{Mahler}. 
As implied by the following theorem due to Palmeira, foliations of the
type $(M,\F) \times \R$ include an important class of plane foliations
(i.e.~foliations whose leaves are diffeomorphic to some Euclidean
space).

\begin{thm}[\cite{CFBP}] 
If $\F$ is a transversely orientable plane foliation on an orientable
$n$-dimensional manifold $M$ (with $n\geq 3$), with finitely generated
fundamental group, such that all leaves are closed, then there exists
a two dimensional surface $S$ and a plane foliation $\F_0$ on $S$ such
that $\F$ is conjugate to the product of $\F_0$ with $\R^{n-2}$.
\end{thm}

Another class of examples is described in the following lemma.

\begin{lem} Let $\pi : M \to B$ be a locally trivial fibration with
compact fiber $L$ and open base, and let $\F$ be a foliation on $M$ that
is transverse to the fibers of $\pi$, in the sense that $T_x\F +
\dim{Ker}\pi_{*_x} = T_xM$ for all $x$ in $M$. Then the \fm $(M,\F)$
is open. 
\end{lem}

\Pf Let $g : B \to [0,\infty)$ be a proper Morse function without local
maxima. It is sufficient to prove that if $f : M \to [0,\infty)$ is
sufficiently close to $\pi^*g$ in the fine $C^\infty$ topology, then
$f$ is proper and has no \lw local maxima. It is easy to see that a
function that is $C^0$-close to a proper function is proper as well. Let 
$\{x^j ; j \geq 1 \}$ be the set of critical points of $g$. For every 
$j$, let $x^j_1,\dots,x^j_k$ be local coordinates about $x^j$, defined 
on a \nd $U^j$ of $x^j$ in $B$, for which  
$$
g(x^j_1, \ldots , x^j_k) = (x^j_1)^2 + g(0, x^j_2, \ldots , x^j_k)\;.
$$
Assume also that $U^j$ is small enough to guarantee existence of a 
trivialization $\phi : \pi ^{-1}(U^j) \to U^j \times L$, with the 
property that the image of each local sections $\sigma^j_\ell : U^j
\to \pi^{-1}(U^j) : x \mapsto \phi^{-1}(x,\ell)$, $\ell \in L$ is
contained in a leaf of $\F$. A function $f$ whose $2$-jet is
sufficiently close to that of $\pi^*g$ satisfies the following
properties.   
\begin{enumerate}
\item[-] The \lw critical points of $f$ are all contained in $\pi^{-1}
(\cup_j U^j)$.   
\item[-] For all $j\geq 1$, and for all $\ell$ in $L$, the second
derivative of $f \circ \sigma^j_\ell$ in the direction of $x^j_1$
is strictly positive on $U^j$.   
\end{enumerate} 
Then, if $y$ is a \lw critical point of $f$, the foliated second
differential of $f$ at $y$ may not be negative definite, that is, $y$
may not be a \lw local maximum of~$f$.  
\cqfd

\begin{rmk}{\rm
It is not know to us whether locally trivial fibrations with open
fibers are always open as foliated manifolds. Nevertheless, such
foliated manifolds satisfy the conclusion of \tref{ht}, as implies by
\tref{Bach} (cf.~\cite{B, B-f}). 
}\end{rmk}

We will now describe a general procedure that allows one to construct
open foliated manifolds. We need to recall the notion of Novikov
component of a codimension one foliation $\F$ on a closed
manifold~$M$. 

\begin{dfn}[\cite{Nov}] Two points $x$ and $y$ in $M$ are said to be
{\em equivalent \wrt $\F$} if either $F_x = F_y$, or the foliation 
$\F$ admits a closed transversal that contains both $x$ and $y$.   
\end{dfn}
It is not difficult to verify that this defines an equivalence
relation on~$M$.  

\begin{dfn} A {\em Novikov component of the foliation $\F$} is an
equivalence class for this equivalence relation.
\end{dfn}

\begin{ex}{\rm The Reeb foliation on $S^3$ has three Novikov
components. One of them is the torus leaf. The other two are two open solid 
tori bounded by the torus leaf.  
}\end{ex}

\begin{thm}[\cite{Nov}] A Novikov component is either a compact leaf or 
an open submanifold whose boundary is a finite union of compact leaves 
which are themselves Novikov components. 
\end{thm}
The following result is due to Ferry and Wasserman.

\begin{thm}[\cite{F-W}] For a codimension one foliation $\F$ on a closed 
manifold $M$, the following statements are equivalent.
\begin{enumerate}
\item[-] $\F$ has one Novikov component.
\item[-] There exists an $\F$-generic function $f : M \to \R$ with no
\lw degenerate critical points.
\end{enumerate}
\end{thm}

Consider now a codimension one foliation $\F$ with one Novikov
component on a closed $(n+1)$-dimensional manifold $M$. Let $f:M\to
\R$ be an $\F$-generic function without \lw degenerate critical
points. Let $\Gamma \subset M$ be the set of leafwise local maxima
of $f$. Notice that since $f$ has no \lw degenerate critical points,
its singular locus (\dref{lwcr}) is necessarily transverse to $\F$
(cf.~\sref{TBs} or \cite{F-W}). It is therefore a finite union of
embedded circles intersecting $\F$ transversely; the set $\Gamma$ is
the union of some of them.  

\begin{lem} The \fm $(M' = M - \Gamma, \F' = \F|_{M - \Gamma})$ is open.
\end{lem}

\Pf We already have a bounded below $\F$-generic function $f|_{M'}$
on $M'$ with no \lw local maxima. The only thing that needs to be done 
is to modify $f|_{M'}$ so as to make it proper. Take a \tnd $U$ of $\Gamma$ 
in $M$. The set $U$ is the image of an embedding $e : E \to M$, defined 
on the total space of a rank-$n$ vector bundle $p: E\to \Gamma$, such that
$e \circ s = \dim{Id}_\Gamma$, where $s : \Gamma \to E$ is the zero 
section. Since $\Gamma$ is transverse to $\F$, we may assume, without loss 
of generality, that the foliation $\F|_U$ corresponds, under the map 
$e$, to the foliation of $E$ by the fibers of $p$. Such a \tnd is called 
hereafter a {\em foliated tubular neighborhood}. Let $g$ be a Riemannian 
metric on the bundle $E$. Let $\theta : \R \to [0,1]$ be a function with 
compact support whose only critical value, aside from $0$, is a global 
maximum achieved at the point $0$. Consider the function    
$$
f' : M' \to \R : x \mapsto \left\{ \begin{array}{ll} 
f(x) + \displaystyle{\frac{\theta(g(e^{-1}x, e^{-1}x))}{g(e^{-1}x,
e^{-1}x)}} & \mbox{for} \;\; x \;\;\dim{in}\;\; U - \Gamma \\ 
f(x)  & \mbox{otherwise}\;.
\end{array}\right.$$ 
The function $f'$ is proper and bounded below. Moreover, provided $U$
is small enough, its \lw critical points are exactly those of $f|_{M'}$. 
Thus $f'$ has no \lw local maxima. 
\cqfd

\begin{rmk}{\rm Making $f|_{M - \Gamma}$ proper requires the set
$\Gamma$ of \lw local maxima of $f$ to be a union of circles. If on
the contrary, $\Gamma$ has a line segment as one of its connected 
components, as might be the case if the function $f$ had leafwise 
degenerate critical points, it would not be possible to make $f|_{M - 
\ol{\Gamma}}$ proper without creating \lw local maxima.   
}\end{rmk}

\begin{rmk}{\rm It is a classical result due to Rummler and Sullivan
that a transversely orientable, codimension one foliation has one
Novikov component \Iff it is geometrically taut, that is, \Iff $M$
admits a Riemannian metric \wrt which all the leaves of $\F$ are
minimal submanifolds (a proof of this result can be found
in~\cite{C-C}).   
}\end{rmk}

As was suggested to us by Takashi Tsuboi, interesting foliations with
one Novikov component are obtained by forming connected sums along
closed trans\-ver\-sals. Let $M_1$ and $M_2$ be $(n+1)$-dimensional
manifolds endowed with codimension one foliations $\F_1$ and $\F_2$
respectively. For $i = 1,2$, let $c_i : S^1 \to M_i$ be an embedding
transverse to $\F_i$. Suppose that the normal bundle of $c_i(S^1)$ is
trivial (which it is when the foliated manifold $(M_i,\F_i)$ is
orientable). Consider a foliated \tnd $U_i$ of $c_i(S^1)$. The set
$U_i$ is the image of an embedding $e_i : S^1 \times \R^n \to M_i$
that coincides with $c_i$ on $S^1 \times \{0\}$, and that maps the
``vertical'' foliation of $S^1 \times \R^n$ isomorphically onto the
foliation $\F_i|_{U_i}$. {\em The connected sum of the foliated
manifolds $(M_1,\F_1)$ and $(M_2,\F_2)$ along the closed transversals
$c_1$ and $c_2$}, denoted by $(M_1,\F_1)_{c_1}\#_{c_2} (M_2,\F_2)$, is
defined to be the \fm obtained by gluing $M_1 - c_1(S^1)$ to $M_2 -
c_2(S^1)$ along $U_1 - c_1(S^1)$ and $U_2 - c_2(S^1)$ via the
isomorphism       
$$
U_1 - c_1(S^1) \to U_2 - c_2(S^1) : e_1(t,x) \mapsto e_2(t,
\fr{x}{\|x\|^2})\;.   
$$
If both $\F_1$ and $\F_2$ have one Novikov component, then their
connected sum has one Novikov component as well. Notice that one could
as well perform the connected sum along disjoint unions of embedded
circles, $c_1^1 \cup \ldots \cup c_1^k$ in $M_1$, and $c_2^1 \cup
\ldots \cup c_2^k$ in $M_2$.\\     

Consider, for instance, a  manifold of the type $M = S^1 \times F$,
where $F$ is some $n$-dimensional manifold. Let $\pi$ denote the natural
projection $M \to S^1$, and let $\F$ denotes the (trivial) foliation
of $M$ by the fibers of $\pi$. Let also $c_1, c'_1, c_2 : S^1 \to 
M$ be disjoint embeddings such that $\pi \circ c_1 = \pi \circ c'_1 = 
\dim{Id}_{S^1}$, and such that $\pi \circ c_2$ coincides with the double 
cover $e^{i t} \mapsto e^{i 2 t}$ (alternatively, one could suppose that 
$\pi \circ c_1 = \pi \circ c'_1$ and $\pi \circ c_2$ are different
covers of $S^1$ of 
the type $e^{i t} \mapsto e^{i n t}$). Then the connected sum of
$(M,\F)$ with itself along the closed transversals $c_1 \cup c'_1$ and
$c_2 \cup c'_1$ is a \fm with one Novikov component, and with dense
leaves. Indeed, a leaf of the connected sum corresponds to an
equivalence class for the equivalence relation on $S^1$ generated by
$e^{i t} \sim e^{i 2 t}$. The open \fm obtained by removing from $M$
the set of closed curves along which some given $\F$-generic function
without any \lw degenerate critical point achieves \lw local
maxima has dense leaves as well. Besides, if $F$ is an open manifold,
it is not necessary to remove something to achieve openness.     

\begin{lem} If the manifold $F$ is open, then the \fm $(M,\F)_{c_1 
\cup c'_1} \\ \#_{c_2 \cup c'_1} (M,\F)$ is open. 
\end{lem} 

\begin{figure}[h]
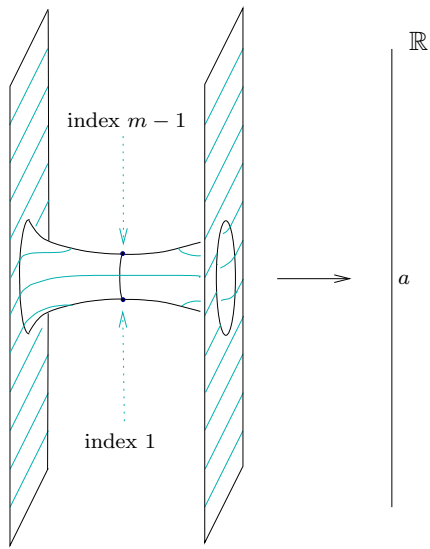

\begin{center}
\input h-p-connected-sum.pstex_t
\smallskip\noindent 
\caption{The function $g$ on the connected sum of two corresponding 
fibers} 
\label{fig:connected-sum}
\end{center}
\end{figure}

\Pf Let $f : F \to [0,\infty)$ be a proper Morse function without local 
maxima. Let $a, b$ be distinct noncritical values of $f$. Let $\gamma :  
S^1 \to F$ be an injective map whose image is contained in the level 
$f^{-1}(a)$. Let $q$ (\rp $q'$) be a point in $f^{-1}(a) - \gamma(S^1)$ 
(\rp $f^{-1}(b)$). Suppose that the closed transversal $c_1$ (\rp $c'_1$) 
coincides with the map $S^1 \to M : t \mapsto (t,q)$ (\rp $S^1 \to M : 
t \mapsto (t,q')$), and that $c_2$ coincides with the map $S^1 \to  
M : t \mapsto (2t,\gamma(t))$. Suppose also that on the \tnds $U_1$, 
$U'_1$ and $U_2$ of $c_1$, $c'_1$ and $c_2$ respectively, along which
the connected sum is performed, the function $p^*f$ corresponds to some
function  
$$
S^1 \times \R^n \to \R : (t,x_1, \dots, x_n) \mapsto c + h(x_n)\;,
$$
where $p$ denotes the projection $M \to F$, where $h : \R \to \R$ is
an embedding, and where $c$ is either $a$ or $b$. The functions $f|_{M
- (c_1(S^1) \cup c'_1(S^1))}$ and $f|_{M - (c_2(S^1) \cup c'_1(S^1))}$
assemble into a proper function $g : M_{c_1 \cup c'_1} \#_{c_2 \cup
c'_1} M \to [0,\infty)$ that coincides with $f$ on $(M - (U_1 \cup
U'_1)) \coprod (M - (U_2 \cup U'_1))$, and whose critical locus in
$U_1 - (c_1(S^1) \cup c'_1(S^1)) \simeq U_2 -  (c_2(S^1) \cup
c'_1(S^1))$ is made of two disjoint transverse closed curves,
consisting of \lw critical points of \lw index $m-1$ and $1$
respectively. Figure~\ref{fig:connected-sum} shows how to construct
$g$ on the connected sum of a fiber in the \tnd $U_1$ (or $U'_1$) with
the corresponding fiber in the \tnd $U_2$ (or $U'_1$).   
\cqfd

\section{Application to Poisson geometry}\label{Poisson}

Let us recall that a regular Poisson manifold can be described as a 
foliated space $(M,\F)$ endowed with a \lw symplectic structure, that
is, a section of the second exterior power of the cotangent bundle
$T^*\F$ of the foliation $\F$, whose restriction to each leaf of $\F$
is a symplectic form (cf.~\cite{AW-book}). The question of existence
of such a structure on a given \fm as been approached in a previous
paper~\cite{B-f} (see also \cite{B}), where examples of foliations are
presented that do not support any \lsst although the obvious
obstructions vanish. On the other hand, as explained below, a \lsst on
the \fm $(M,\F)$ is a solution of a certain open, foliated invariant
differential relation.\\      

Consider a \fm $(M,\F)$. A section of the bundle $\Lambda^kT^*\F$ is 
called a {\em \td $k$-form}. A \td $2$-form $\alpha$ is said to be
{\em nondegenerate} if for every point $x$ in $M$, the skewsymmetric
bilinear form $\alpha(x)$ on $T_x\F$ is nondegenerate. The usual
exterior differential restricts naturally into a map $d_\F :
\Gamma(\Lambda^\star T^*\F) \to \Gamma( \Lambda^{\star +1} T^*\F)$
(cf.~\cite{HMS}). The corresponding cohomology  
$$
H^{\star}(\F) = \displaystyle{\frac{\{d_\F-\dim{closed} 
\star-\dim{forms}\}}{\{d_\F-\dim{exact} \star-\dim{forms}\}}}\;,
$$
is called the {\em tangential de~Rham cohomology of the \fm
$(M,\F)$}. There is a natural affine fibration $L : (T^*\F)^1 \to
\Lambda^2 T^*\F $ defined by $L(j^1{\alpha}(x)) = {d_\F\alpha}(x)$. 
Choose a $d_\F$-closed \td $2$-form ${\theta}$. Define $\Omega_\theta$
to be the differential relation $\Omega_{\theta} = \{\, j^1{\alpha}(x)
\,; \,{\theta}(x) + L(j^1{\alpha}(x))$ is nondegenerate $\}$. A
solution of $\Omega_{{\theta}}$ is thus a tangential differential
$1$-form $\alpha$ such that ${d_\F\alpha} + {\theta}$ is a \lsst
belonging to the same tangential de~Rham cohomology class
as~${\theta}$. Thus, the set of solutions, modulo the set of
$d_\F$-closed $1$-forms, parameterizes the set of \lsste s lying in the
class $[\theta]$. On the other hand, since the map $L$ is an affine
fibration, the space of sections of $\Omega_\theta$ is weakly homotopy
equivalent to the space of \lw nondegenerate $2$-forms.

\begin{prop}[\cite{Giroux}]
The relation $\Omega_{\theta}$ is open and invariant under foliated isotopies. 
\end{prop}

\Pf That $\Omega_\theta$ is open follows directly from its
definition. The second assertion relies on considering the ``right''
lift for isotopies, that~is~:  
$$
\varphi_{t}\cdot\alpha = \varphi_{t}^{*}\alpha + 
\int^t_0\varphi^*_s\left(i(X_{s})\theta\right)\, ds\;,
$$
where $\alpha$ is an element of $T^{*}\F$, and where $X_{t}$ denotes
the time-dependent vector field associated to the isotopy
$\varphi_{t}$. If $\alpha$ is a tangential differential $1$-form for
which the tangential differential $2$-form ${d_\F\alpha} + \theta$ is
\lw nondegenerate, then $\varphi^*_t({d_\F\alpha} + \theta)$ is \lw
nondegenerate as well (while the form $\varphi^*_t (d_\F \alpha) +
\theta$ might very well be leafwise degenerate),~and     
\begin{eqnarray*}
\varphi^*_t({d_\F\alpha} + \theta) & = &  d_\F\varphi_{t}^{*}\alpha +
\theta + (\varphi^*_t\theta - \theta) \\
{} & = &  d_\F\varphi^*_t\alpha + \theta + 
\int^t_0\left.\frac{d}{ds}\varphi^*_s\theta\right|_{s} \, ds \\
{} & = &  d_\F\varphi^*_t\alpha + \theta + \int^t_0\varphi ^*_s 
\EK_{X_{s}}\theta\, ds \\
{} & = &  d_\F\varphi^*_t\alpha + \theta + \int^t_0
\varphi^*_s  d_\F i(X_{s})\theta\, ds \\
{} & = &  d_\F\left( \varphi^*_t\alpha 
+ \int^t_0 \varphi^*_s i(X_{s})\theta\, ds  \right) + 
\theta\; .
\end{eqnarray*}
\cqfd

\vspace{.2cm}
Thus, \tref{ht} applies to the relation $\Omega_\theta$.
 
\begin{thm}\label{Richter}
Let $(M,\F)$ be an open foliated manifold. Given a $d_\F$-closed \td
$2$-form $\theta$, any family $\beta_s, s\in [0,1]^p$ of \lw nondegenerate 
$2$-forms is homotopic to a family of \lsste s cohomologous to $\theta$. 
Moreover, if $\beta_s$ is already symplectic and cohomologous to $\theta$ 
for $s$ in $\partial ([0,1]^p)$, the homotopy may be chosen to be stationary 
for those parameters. Equivalently, the inclusion of the space of \lsste s 
cohomologous to $\theta$ into the space of \lw nondegenerate $2$-forms 
is a weak homotopy equivalence 
\end{thm}

\tref{Richter} implies the following existence and uniqueness result 
for \lw \s structures.

\begin{cor}
On an open foliated manifold any \lw nondegenerate $2$-form can be
deformed into a \lw symplectic form (with prescribed \tdco
class). Moreover, if two cohomologous \lsste s can be joined by a path
of \lw nondegenerate $2$-forms, they can also be joined by a path of
cohomologous \lw symplectic forms.    
\end{cor}


\begin{thebibliography}{xx}

\bibitem{B} M. Bertelson, {\em Foliations associated to regular 
Poisson structures}, Ph.D. Thesis, Stanford University, June 2000. 

\bibitem{B-f} M. Bertelson, {\em Foliations associated to regular 
Poisson structures}. Accepted for publication in
Commun. Contemp. Math. 

\bibitem{Boardman}J.M. Boardman, Singularities of differentiable
maps. {\em Inst. Hautes \'Etudes Sci. Publ. Math. No.} {\bf 33} (1967),
21--57.

\bibitem{C-C} A. Candel and L. Conlon, {\em Foliations. I.}
Graduate Studies in Mathematics, 23. American Mathematical Society,
Providence, RI, 2000.

\bibitem{AW-book} A. Cannas da Silva and A. Weinstein, {\em Geometric
Models for Noncommutative Algebras.} Berkeley Mathematics Lecture
Notes Series. American Mathematical Society, 1999.

\bibitem{F-W} S. Ferry and A. Wasserman, Morse theory for
codimension-one foliations. {\em Trans. Amer. Math. Soc.} 298 (1986),
no. 1, 227--240.  

\bibitem{Giroux} E. Giroux, {\em Flexibilit\'e en g\'eom\'etrie
symplectique, d'apr\`es Mikhael Gromov.} Unpublished lecture notes
(1993). 

\bibitem{GG} M. Golubitsky and V. Guillemin, {\em Stable mappings and
their singularities.} Graduate Texts in Mathematics,
Vol. 14. Springer-Verlag, New York-Heidelberg, 1973. 

\bibitem{Gv-69} M. Gromov, Stable mappings of foliations into
manifolds. {\em Izv. Akad. Nauk SSSR Ser. Mat.} {\bf 33} (1969),
707--734. 

\bibitem{Gv-book} M. Gromov, {\em Partial differential
relations.} Ergebnisse der Mathematik und ihrer Grenzgebiete (3)
[Results in Mathematics and Related Areas (3)], 9. Springer-Verlag,
Berlin-New York, 1986. 

\bibitem{HMS} G. Hector, E. Mac\'\i as, M. Saralegi, Lemme de Moser
feuillet\'e et classification des vari\'et\'es de Poisson
r\'eguli\`eres. (French) [A foliated version of Moser's lemma and
classification of regular Poisson manifolds] {\em Publ. Mat.} {\bf 33}
(1989), no. 3, 423--430.  

\bibitem{Igusa} K. Igusa, Higher singularities of smooth functions 
are unnecessary. {\em Ann. of Math. (2)}, {\bf 119} (1984), 1--58.

\bibitem{JMt} J. Martinet, {\em Singularities of smooth functions 
and maps.} Translated from the French by Carl P. Simon. London 
Mathematical Society Lecture Note Series, 58. 
Cambridge University Press, Cambridge-New York, 1982. 

\bibitem{JM-73} J.N. Mather, On Thom-Boardman singularities.
Dynamical systems (Proc. Sympos., Univ. Bahia, Salvador, 1971),
233--248. Academic Press, New York, 1973. 

\bibitem{JM-68} J.N. Mather, Stability of $C\sp{\infty }$ mappings. III. 
Finitely determined mapgerms. Inst. Hautes Études Sci. Publ. Math. No. 35, 
1968 279--308. 

\bibitem{Milnor} J. Milnor, {\em Morse theory. Based on lecture notes
by M. Spivak and R. Wells.} Annals of Mathematics Studies, No. 51
Princeton University Press, Princeton, N.J. 1963. 

\bibitem{Nov} S.P. Novikov, Topology of foliations. {\em Trans. Moscow
Math. Soc.} {\bf 14} (1965), 268--304.

\bibitem{CFBP} C.F.B. Palmeira, Open manifolds foliated by
planes. {\em Ann. Math. (2)}  {\bf 107} (1978), no. 1, 109--131.

\bibitem{CT} C.H. Taubes, The Seiberg-Witten invariants and symplectic
forms. {\em Math. Res. Lett.} {\bf 1} (1994), no. 6, 809--822.

\bibitem{Whitney} H. Whitney, Elementary structure of real 
algebraic varieties. {\it Ann. of Math. (2)} {\bf 66} (1957), 545--556.

\end{thebibliography}
\end{document}